\documentclass[review]{elsarticle}

\usepackage{lineno,hyperref}
\usepackage{a4wide}
\usepackage{amsmath}
\usepackage{amsthm}
\usepackage{amssymb}
\usepackage{amsfonts}
\usepackage{mathrsfs}
\usepackage{stmaryrd}
\usepackage{graphicx}
\usepackage{epstopdf}
\usepackage{paralist}
\usepackage{subfigure}
\usepackage{xcolor}
\usepackage{enumitem}
\usepackage{units}
\newcommand{\R}{\mathbb{R}}

\newcommand{\calD}{\mathcal{D}}
\newcommand{\calE}{\mathcal{E}}

\newcommand{\calP}{\mathcal{P}}

\newcommand{\calT}{\mathcal{T}}
\newcommand{\calV}{\mathcal{V}}

\newcommand{\n}{\mathbf{n}}

\newcommand{\meas}[1]{|{#1}|}

\newcommand{\ff}{\mathrm{ff}}
\newcommand{\porm}{\mathrm{pm}}
\newcommand{\ipmff}{\mathrm{if}}
\newcommand{\vel}{\mathbf{v}}
\newcommand{\BFCoeff}{\alpha_{\mathrm{BF}}}
\newcommand{\vertex}{v}
\newcommand{\velx}{v_x}
\newcommand{\vely}{v_y}
\newcommand{\veli}[1]{v_{#1}}
\newcommand{\up}{\mathrm{up}}

\theoremstyle{plain}

\newtheorem{definition}{Definition}

\theoremstyle{remark}

\newcommand{\Dumux}{{Du\-Mu$^\text{x}$ }}
\newcommand{\eqbydef}{:=}

\newcommand*\patchAmsMathEnvironmentForLineno[1]{%
  \expandafter\let\csname old#1\expandafter\endcsname\csname #1\endcsname
  \expandafter\let\csname oldend#1\expandafter\endcsname\csname end#1\endcsname
  \renewenvironment{#1}%
  {\linenomath\csname old#1\endcsname}%
  {\csname oldend#1\endcsname\endlinenomath}}%
\newcommand*\patchBothAmsMathEnvironmentsForLineno[1]{%
  \patchAmsMathEnvironmentForLineno{#1}%
  \patchAmsMathEnvironmentForLineno{#1*}}%
\AtBeginDocument{%
  \patchBothAmsMathEnvironmentsForLineno{equation}%
  \patchBothAmsMathEnvironmentsForLineno{align}%
  \patchBothAmsMathEnvironmentsForLineno{flalign}%
  \patchBothAmsMathEnvironmentsForLineno{alignat}%
  \patchBothAmsMathEnvironmentsForLineno{gather}%
  \patchBothAmsMathEnvironmentsForLineno{multline}%
}

\journal{Journal of Computational Physics}

\begin{document}

\begin{frontmatter}

\title{Coupling Staggered-Grid and MPFA Finite Volume Methods for Free Flow/Porous-Medium Flow Problems}


\author[addressIWS]{Martin Schneider\corref{mycorrespondingauthor}}
\ead{martin.schneider@iws.uni-stuttgart.de}
\author[addressIWS]{Kilian Weishaupt}
\ead{kilian.weishaupt@iws.uni-stuttgart.de}
\author[addressIWS]{Dennis Gl\"aser}
\ead{dennis.glaeser@iws.uni-stuttgart.de}
\author[addressIWS]{Wietse M. Boon}
\ead{wietse.boon@iws.uni-stuttgart.de}
\author[addressIWS]{Rainer Helmig}
\ead{rainer.helmig@iws.uni-stuttgart.de}

\cortext[mycorrespondingauthor]{Corresponding author}

\address[addressIWS]{Institute for Modelling Hydraulic and Environmental Systems,
           University of Stuttgart,
             Pfaffenwaldring 61,
               70569 Stuttgart, Germany}

\begin{abstract}
  A discretization is proposed for models coupling free flow with anisotropic porous medium flow. Our approach
  employs a staggered grid finite volume method for the Navier-Stokes equations in the free flow
  subdomain and a MPFA finite volume method to solve Darcy flow in the porous medium. After
  appropriate spatial refinement in the free flow domain, the degrees of freedom are conveniently
  located to allow for a natural coupling of the two discretization schemes. In turn, we
  automatically obtain a more accurate description of the flow field surrounding the porous medium.
  Numerical experiments highlight the stability and applicability of the scheme in the presence of
  anisotropy and show good agreement with existing methods, verifying our approach.
\end{abstract}


\end{frontmatter}


\section{Introduction}
\label{sec:introduction}
Coupled systems of free flow and flow through a porous medium can be found ubiquitously in various kinds of natural and industrial contexts, including
soil water evaporation \cite{vanderborght2017a}, fuel cell water management \cite{gurau2009a}, food processing \cite{verboven2006a}, and evaporative cooling for turbomachinery \cite{dahmen2014a}.
Despite every-growing computational capacities, a discrete numerical modeling these kind of systems, including the porous geometry, is only feasible for small-scale problems. For larger domain
sizes, averaging techniques involving the concept of a representative elementary volume (REV) \cite{bear1972a}, are used in order to yield upscaled models for the description of the porous
domain \cite{whitaker1999a}. These models can then be coupled to the free-flow region either using a single-domain or a two-domain approach. For the former, both the porous medium and the
free flow are described by a single set of equations, as first introduced by Brinkman \cite{brinkman1949a}. The two domains are then discerned by a spatial variation of material parameters.
On the other hand, the two-domain approach decomposes the problem into two disjoint subdomains. The free-flow region is then governed by the Navier-Stokes equations while Darcy's or Forchheimer's law is used in the porous medium subdomain
\cite{ochoa1995a, layton2002, jamet2009a, mosthaf2011a}. In order to maintain thermodynamic consistency, appropriate coupling conditions have to be formulated which enforce the conservation
of mass, momentum and energy across the interface between the two domains \cite{hassanizadeh1989a}. 
The aim of this work is to couple two different discretization schemes, and we will focus on the two-domain approach since it provides this flexibility more readily.
We will consider laminar single-phase flow in the following but the extension to compositional multi-phase flow in the porous domain \cite{mosthaf2011a} and the use of turbulence models in the free flow part \cite{fetzer2017a} is possible.

Being an active area of research, various mathematical and numerical models for the
coupling of (Navier-) Stokes and Darcy/Forchheimer have been developed during recent years. 
Examples using the same spatial discretization scheme for both domains range from the staggered-grid finite volume method \cite{iliev2004a, harlow1965a},
finite element method \cite{discacciati2009a} or the box-scheme, a vertex-centered finite volume method \cite{mosthaf2011a, huber2000a}. Furthermore,
combinations of colocated and staggered-grid schemes have been developed \cite{rybak2015a, masson2016a, fetzer2017a}. For more discretization approaches
 to the modeling of coupled free flow and porous medium flow, we refer the reader to \cite{discacciati2009a} and references therein. 

In this work, the free flow equations are discretized using the staggered-grid method, which forms a stable numerical method for such problems. In turn, no additional stabilization techniques are necessary, in contrast to colocated schemes \cite{versteeg2007a}, for example. In the porous medium, we employ a multi-point-flux-approximation method (MPFA) \cite{aavatsmark2002a}. This method was developed to overcome the shortcomings of the classical two-point-flux approximation.
In particular, the MPFA scheme does not require the grid to be $\mathbf{K}$-orthogonal, meaning that the grid cells need not be in line with the principal directions of the permeability
tensor $\mathbf{K}$. This is especially important for skewed and unstructured grids or in case the principal directions of the permeability tensor are inclined, in the case of
layered porous structures or faults, for example. The MPFA method has been applied previously for solving Brinkman's equation \cite{iliev2014a, brinkman1949a} for a coupled system of free flow and porous medium. The referenced works thus use the same discretization scheme for both subdomains.
The novelty of this work is that we employ a staggered-grid method for the free-flow model (Navier-Stokes), and couple it with a MPFA method for the porous medium model (Darcy). 

The paper is organized as follows. Section~\ref{sec:equations} introduces the coupled model by presenting the governing equations in the free flow and porous medium flow subdomains, respectively, as well as the coupling conditions at the interface. 
The discretization schemes are introduced in Section~\ref{sec:discretization}, including the newly proposed coupling at the interface between the two subdomains. 
Numerical results are presented in Section~\ref{sec:numresults} with the use of two test cases. 
Finally, Section~\ref{sec:conclusions} focuses on the conclusions.

\section{Governing Equations}
\label{sec:equations}

In order to present the governing equations for the coupled model, we first introduce the assumptions on the geometry and the notational conventions. 
After these preliminary definitions, we continue with the equations governing free flow and porous medium flow, followed by the coupling conditions.

Let $\Omega\subset\R^d$, $d\in \{2, 3\}$, be an open, connected
Lipschitz domain with boundary $\partial\Omega$ and $d$-dimensional measure $\meas{\Omega}$.
Furthermore, let $\Omega^\porm$ and $\Omega^\ff$ be a disjoint partition of $\Omega$ representing the porous medium and free flow subdomains, respectively. 
The subdomain boundaries are given by the interface $\Gamma^\ipmff \eqbydef \partial \Omega^\ff \cap \partial \Omega^\porm$ as well as the remainders $\Gamma^\porm \eqbydef \partial \Omega^\porm \setminus \Gamma^\ipmff$ and $\Gamma^\ff \eqbydef \partial \Omega^\ff \setminus \Gamma^\ipmff$. For brevity, the superscripts are often omitted when no ambiguity arises. 

The external boundary $\partial \Omega = \Gamma^\ff \cup \Gamma^\porm$ is further decomposed such that $\Gamma^\ff = \Gamma_\mathrm{v}^\ff \cup \Gamma_\mathrm{p}^\ff$ and $\Gamma^\porm = \Gamma_\mathrm{v}^\porm \cup \Gamma_\mathrm{p}^\porm$ disjointly. Here, the subscript denotes whether the velocity or the pressure is prescribed as a boundary condition. 
To ensure unique solvability of the resulting system, we assume that $\meas{\Gamma_\mathrm{p}^\porm \cup \Gamma_\mathrm{p}^\ff} > 0$, i.e. that a pressure boundary condition is imposed on a subset of the boundary $\partial \Omega$ with positive measure.

Let $\n$ denote the unit normal vector on $\partial \Omega$ oriented outward with respect to $\Omega$. We abuse notation and let $\n$ moreover denote the unit normal vector on $\Gamma^\ipmff$ oriented outward with respect to $\Omega^\ff$.

\subsection{Free Flow}
\label{sub:free_flow}

In our model, the Navier-Stokes equations govern the free flow in subdomain $\Omega^\ff$. These equations are given by:
\begin{subequations}
  \label{eq:navierstokes}
  \begin{align}
      \frac{\partial \varrho}{\partial t} + \nabla \cdot (\varrho \vel) &= q, \\
  \frac{\partial (\varrho \vel)}{\partial t} + \nabla \cdot \left(\varrho \vel \vel^{\mathrm{T}}
      - \mu (\nabla \vel + (\nabla \vel)^{\mathrm{T}})
      + p \mathbf{I} \right)
      &= \varrho \textbf{g}, &&\mathrm{in} \, \Omega^\ff, \\
       \vel &= \vel_\Gamma, &&\mathrm{on} \, \Gamma^\ff_\mathrm{v}, \\
       p &= p_\Gamma, &&\mathrm{on} \, \Gamma^\ff_\mathrm{p}.
   \end{align}
\end{subequations}

The unknown variables are the velocity $\vel$ and the pressure $p$. Here, $\varrho$ and $\mu$ denote the potentially pressure-dependent density and viscosity, respectively, while $q$ is a source (or sink) term, and $\mathbf{g}$ describes the influence of gravity. $\mathbf{I}$ is the identity tensor in $\R^{d \times d}$. The boundary conditions are given by known quantities $\vel_\Gamma$ and $p_\Gamma$, representing the velocity or pressure at the corresponding boundary.

\subsection{Porous Medium Flow}
\label{sub:porous_medium_flow}

The equations governing single-phase flow in the porous-medium $\Omega^\porm$ are given by
\begin{subequations}
    \label{eq:darcy}
  \begin{align}
      \frac{\partial \varrho}{\partial t} + \nabla \cdot \left( \varrho \vel \right) &= q, \\
      \vel + \frac{1}{\mu} \mathbf{K} \left( \nabla p - \varrho \mathbf{g}\right) &= 0, &&\mathrm{in} \, \Omega^\porm, \label{eq: darcys law}\\
   \vel \cdot \mathbf{n} &= v_\Gamma, &&\mathrm{on} \, \Gamma^\porm_\mathrm{v}, \\
       p &= p_\Gamma, &&\mathrm{on} \, \Gamma^\porm_\mathrm{p}.
  \end{align}
\end{subequations}
Equation \eqref{eq: darcys law} states that the momentum balance in the porous medium is given by Darcy's law, i.e. the Darcy velocity is calculated as $\vel = -\frac{1}{\mu} \mathbf{K} \left( \nabla p - \varrho \mathbf{g}\right)$, with $\mathbf{K}$ being the permeability tensor. Similar to  the free-flow equations \eqref{eq:navierstokes}, $q$ denotes a source or sink term. Finally $v_\Gamma$ and $p_\Gamma$ are known quantities representing the normal flux and pressure on the corresponding boundaries, respectively.

\subsection{Coupling Conditions}
\label{sub:coupling}

In order to derive a thermodynamically consistent formulation of the coupled problem,
 conservation of mass and momentum has to be guaranteed at the interface between the porous
 medium and the free-flow domain. We therefore impose the following interface conditions:

\begin{subequations}\label{eq:interfaceConditions}
 \begin{align}
        \label{eq:interfaceConditionsFlux}
    \vel^\ff \cdot \n &=  - \vel^\porm \cdot \n , \\
    \label{eq:interfaceConditionsMomentum}
    \left(\varrho \vel \vel^\mathrm{T} - \mu (\nabla \vel + (\nabla \vel)^\mathrm{T}) + p \mathbf{I} \right)^\ff \n &= 
      p^\porm \n, \\
     \label{eq:interfaceConditionsBJS}
     \left( - \frac{\sqrt{\mathbf{t} \cdot \mathbf{K} \mathbf{t}}}{\BFCoeff} (\nabla \vel) \mathbf{n} - \vel \right)^\ff \cdot \mathbf{t} &=  0, && \text{on} \, \Gamma^\ipmff. 
 \end{align}
\end{subequations}

The momentum transfer normal to the interface is given by \eqref{eq:interfaceConditionsMomentum} \cite{layton2002}.
Condition \eqref{eq:interfaceConditionsBJS} is the commonly used Beavers-Joseph-Saffman slip condition \cite{beavers1967a, saffman1971a}. Here, $\mathbf{t}$ denotes any unit vector from the tangent bundle of $\Gamma^\ipmff$ and $\BFCoeff$ is a parameter to be determined experimentally. We remark that this condition is technically a boundary condition for the free flow, not a coupling condition between the two flow regimes. Furthermore, it has been developed for free flow strictly parallel to the
interface and might lose its validity for other flow configurations. 

\section{Discretization}
\label{sec:discretization}

This section is devoted to giving an outline of the numerical schemes used in the individual subdomains
 and the incorporation of the interface conditions \eqref{eq:interfaceConditions}. However, we first introduce some notational conventions concerning the partition of $\Omega$ in the following definition.

\begin{definition}[Grid discretization]
\label{def:griddisc}
The tuple \mbox{$\calD \eqbydef (\calT,\calE,\calP,\calV)$} denotes the grid discretization, in which
  \begin{enumerate}[label=(\roman*)]
  \item $\calT$ is the set of control volumes (cells) such that
    \mbox{$\overline{\Omega}= \cup_{K \in \calT} \overline{K}$}.
    For each cell $K\in\calT$, $\meas{K}>0$ denotes the cell volume.
 \item
    $\calE$ is the set of faces such that each face $\sigma$ is a $(d-1)$-dimensional hyperplane
    with measure
    \mbox{$\meas{\sigma}>0$}.
    For each cell $K \in \calT$, $\calE_K$ is the subset of
    $\calE$ such that \mbox{$\partial K  = \cup_{\sigma \in\calE_K}{\sigma}$}. Furthermore, $\mathbf{x}_\sigma$ denotes
    the face evaluation points and $\n_{K,\sigma}$ the unit vector that is normal to
    $\sigma$ and outward to $K$.

  \item
    $\calP \eqbydef \lbrace \mathbf{x}_K\rbrace_{K \in \calT}$ is the set of \emph{cell centers} (not required to be the barycenters) such that
    $\mathbf{x}_K\in K$ and $K$ is star-shaped with respect to $\mathbf{x}_K$.
    For all $K\in\calT$ and $\sigma\in\calE_K$, let $d_{K,\sigma}$ denote the Euclidean distance between $\mathbf{x}_K$ and $\sigma$.
  
  \item  
    $\calV$ is the set of vertices of the grid, corresponding to the corners of the cells. 
  \end{enumerate}
\end{definition}  

For ease of exposition, we assume $d = 2$, however, this is not a limitation and the model can be readily extended for three dimensions.

\subsection{Staggered grid scheme}
\label{sub:staggered}

A staggered-grid finite volume scheme, also known as MAC scheme \cite{harlow1965a}, is used in the free-flow subdomain.
Here, scalar quantities including pressure and density are stored on the cell centers $\calP$ while the velocity degrees of freedom are located on the primary control volumes' faces $\calE$. 
The resulting scheme is stable, hence oscillation-free solutions are guaranteed without the need for additional stabilization techniques. This is in contrast with colocated schemes, in which all unknowns are defined at the same location \cite{versteeg2007a}. $\calD$ is a uniform Cartesian grid with mesh size $h$ in both directions.

\begin{figure} [ht!]
  \centering
  \includegraphics[width=1.0\linewidth,keepaspectratio]{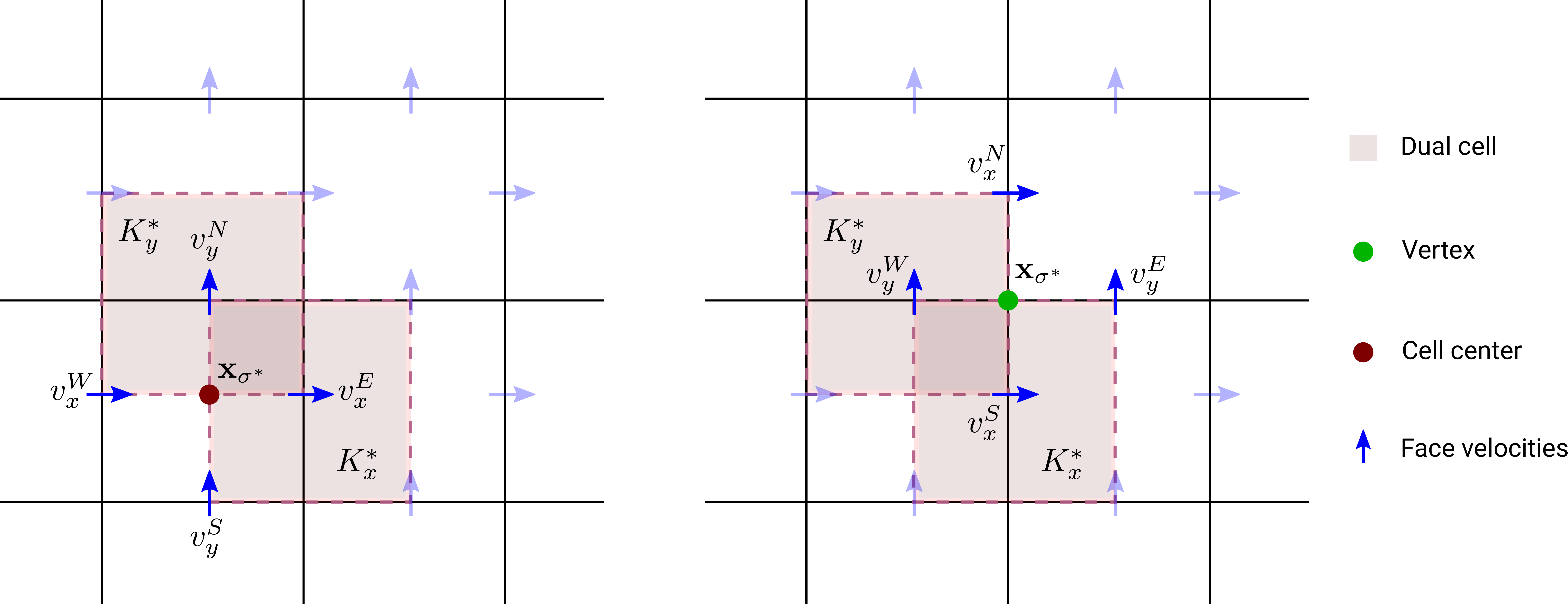}
  \caption{Grid and notations used for the staggered-grid discretization. $K_x^*, K_y^*$ denote the dual cells. The picture on the left illustrates the situation where $\mathbf{x}_{\sigma^*}$ coincides with a cell center, whereas the picture on the right shows the case where the center $\mathbf{x}_{\sigma^*}$ of a dual face $\sigma^*$ coincides with a vertex of the primary grid.}
  \label{pc:staggeredgrid}
\end{figure}

We start with the momentum balance equations. For each face in $\calE$, we construct a secondary control volume $K^*$ with boundary $\partial K^*$, as depicted in Figure~\ref{pc:staggeredgrid}. On each $K^*$, we integrate the momentum balance equation and apply Gauss divergence theorem to obtain
\begin{equation}
  \begin{aligned}
    \int_{K^*} \frac{\partial (\varrho \vel)}{\partial t} \, \mathrm{d} x + \int_{\partial K^*} (\varrho \vel \vel^{\mathrm{T}}) \cdot \n \, \mathrm{d} \Gamma
    &- \int_{\partial K^*} (\mu (\nabla \vel + (\nabla \vel)^{\mathrm{T}})) \cdot \n \, \mathrm{d} \Gamma 
    + \int_{\partial K^*} p \n \, \mathrm{d} \Gamma  
    = \int_{K^*} \varrho \mathbf{g} \, \mathrm{d} x.
  \end{aligned}
\end{equation}

The first and second components of this vector equation are considered separately. Due to the different locations at which the degrees of freedom are defined, we require interpolation operators in order to continue.
Let us introduce the average ($\{ \cdot \}$) and jump quantities ($\llbracket \cdot \rrbracket$) on cell centers ($\calP$) and vertices ($\calV$) of the primal grid for $\vel = [\velx, \vely]^{\text{T}}$ (see Figure \ref{pc:staggeredgrid}):
\begin{equation}
\begin{aligned}
  \{ \vel \}|_{\calP} &= \frac{1}{2} 
  \begin{bmatrix}
    \velx^E + \velx^W \\ \vely^N + \vely^S
  \end{bmatrix}
  & \quad
  \{ \vel \}|_{\calV} &= \frac{1}{2} 
  \begin{bmatrix}
    \velx^N + \velx^S \\ \vely^E + \vely^W
  \end{bmatrix} \\
  \llbracket \vel \rrbracket|_{\calP} &= \frac{2}{h} 
  \begin{bmatrix}
    \velx^E - \velx^W \\ \vely^N - \vely^S
  \end{bmatrix}
  & \quad
  \llbracket \vel \rrbracket|_{\calV} &= \frac{1}{h} 
  \begin{bmatrix}
    \velx^N - \velx^S + \vely^E - \vely^W \\ 
    \velx^N - \velx^S + \vely^E - \vely^W
  \end{bmatrix}
\end{aligned}
\end{equation}
Here, the superscript $\{E, N, W, S\}$ refers to the closest degree of freedom East, North, West, or South of the evaluation point, see Figure \ref{pc:staggeredgrid}.
In the following, the superscript ``up'' denotes the upwind quantity relative to the velocity $\vel$.
Moreover, we introduce $\mu^\mathrm{avg}$ such that in all cell centers $\calP$, it denotes the corresponding viscosity whereas in each vertex $\calV$, it is the viscosity averaged over the adjacent cells. 

With the operators defined, we discretize the momentum balance equation for each secondary control volume $K_i^*$, in which the subscript $i \in \{ x, y \}$ denotes whether the control volume surrounds a vertical or horizontal face, respectively. The discretized equation for component $i$ is then given by
\begin{equation}
  \begin{aligned}
    \int_{K_i^*} \frac{\partial (\varrho v_i)}{\partial t} \, \mathrm{d} x 
    + \int_{\partial K_i^*} (\{ \vel \} (\varrho v_i)^\up) \cdot \n \, \mathrm{d} \Gamma
    - \int_{\partial K_i^*} \mu^\mathrm{avg} \llbracket \vel \rrbracket \cdot \n \, \mathrm{d} \Gamma 
    &+ \int_{\partial K_i^*} p n_i \, \mathrm{d} \Gamma
    = \int_{K^*} \varrho g_i \, \mathrm{d} x.
  \end{aligned}
\end{equation}

We emphasize that the boundary integrals are computed numerically using the following rules for a scalar-valued function $f$ and a vector-valued function $\mathbf{f}$:
\begin{subequations}
  \begin{align}
    \int_{\partial K^*} f \, \mathrm{d} \Gamma 
    &= \sum_{\sigma^* \in \mathcal{E}_{K^*}} \meas{\sigma^*} f(\mathbf{x}_{\sigma^*})
    = h \sum_{\sigma^* \in \mathcal{E}_{K^*}} f(\mathbf{x}_{\sigma^*}), \\
    \int_{\partial K^*} \mathbf{f} \cdot \n \, \mathrm{d} \Gamma 
    &= \sum_{\sigma^* \in \mathcal{E}_{K^*}}   \meas{\sigma^*}\mathbf{f}(\mathbf{x}_{\sigma^*}) \cdot \n_{K^*, \sigma^*}
    = h \sum_{\sigma^* \in \mathcal{E}_{K^*}} \mathbf{f}(\mathbf{x}_{\sigma^*}) \cdot \n_{K^*, \sigma^*}.
  \end{align}
\end{subequations}
By definition, each $\mathbf{x}_{\sigma^*}$ will either be a cell center ($\calP$) or a vertex ($\calV$) of the grid. 

Finally, the mass balance is evaluated on each cell of the grid, i.e. we compute for each $K \in \calT$:
\begin{equation}
  \int_{K} \frac{\partial \varrho}{\partial t} \, \mathrm{d} x 
  + \int_{\partial K} \varrho^\up \vel \cdot \n \, \mathrm{d} \Gamma 
  = 
  \int_{K} q \, \mathrm{d} x.
\end{equation}

The above equations fully define the staggered-grid discretization scheme for all internal faces of the grid $\calD$. For incorporation of boundary conditions, we refer the reader to \cite{versteeg2007a}.

\subsection{Cell-centered finite volume scheme}
\label{sub:ccfvm}

A cell-centered finite volume scheme is employed in the Darcy subdomain, i.e. the
 grid elements are used as control-volumes and the degrees of freedom are associated
 with the cell centers. Typically,
 the finite volume formulation is obtained by integrating the first equation of
 \eqref{eq:darcy} over a control volume $K \in \calT$ and by applying the Gauss divergence theorem:
 \begin{equation}
 \int_K \frac{\partial \varrho}{\partial t} \, \mathrm{d}x +   \sum_{\sigma \in \calE_K }\int_{\sigma} \frac{\varrho^\up}{\mu^\up}\left( - \mathbf{K} \left( \nabla p - \varrho \mathbf{g}\right) \right) \cdot \mathbf{n} \, \mathrm{d} \Gamma = \int_K q \, \mathrm{d}x.
     \label{eq:darcyIntegrated}
 \end{equation}

Replacing, the exact fluxes by an  approximation, i.e. \mbox{$F_{K, \sigma} \approx \int_{\sigma} \left( - \mathbf{K}_K \left( \nabla p - \varrho \mathbf{g}\right) \right) \cdot \mathbf{n} \, \mathrm{d} \Gamma$}
 (here $\mathbf{K}_K$ is the value of $\mathbf{K}$ associated with control volume $K$), yields
\begin{equation}
 \int_K \frac{\partial \varrho}{\partial t} \, \mathrm{d}x + \sum_{\sigma \in \calE_K} \frac{\varrho^\up}{\mu^\up} F_{K, \sigma} = Q_K, \quad \forall \, {K \in \calT},
  \label{eq:darcyCCdiscrete}
\end{equation}
where $F_{K, \sigma}$ is the discrete flux through face $\sigma$ flowing out of cell
 $K$, $Q_K \eqbydef \int_K q \, \mathrm{d}x$ is the integrated source/sink term, and $(\cdot)^\up$ denotes upwinding with respect to the sign of the flux $F_{K,\sigma}$. 
 
Finite volume schemes primarily differ in the approximation of the term
$(\mathbf{\mathbf{K}}_K \nabla p) \cdot \mathbf{n}$
(i.e. the choice of the fluxes $F_{K, \sigma}$). The widely used linear two-point flux approximation (TPFA), for example, is a simple but robust scheme.
However, it is well-known that it is inconsistent on grids that are not $\mathbf{K}$-orthogonal
(see e.g.\ \cite{edwards1998finite}). In this work we consider anisotropic permeability tensors in
the porous medium and $\mathbf{K}$-orthogonality of the grid can thus not be guaranteed.
Therefore, we employ a multi-point flux approximation (MPFA) scheme for the formulation of the discrete
fluxes, which has been presented in \cite{aavatsmark2002a}. This particular scheme is termed MPFA-O and
 is only one among many methods that fall into the family of MPFA schemes (\cite{aavatsmark2008compact, edwards2010mpfafps}).
 Please note that we will omit the suffix ``-O'' throughout this document wherever it would affect the readability.
 
 \begin{figure} [ht!]
  \centering
  \includegraphics[width=0.8\linewidth,keepaspectratio]{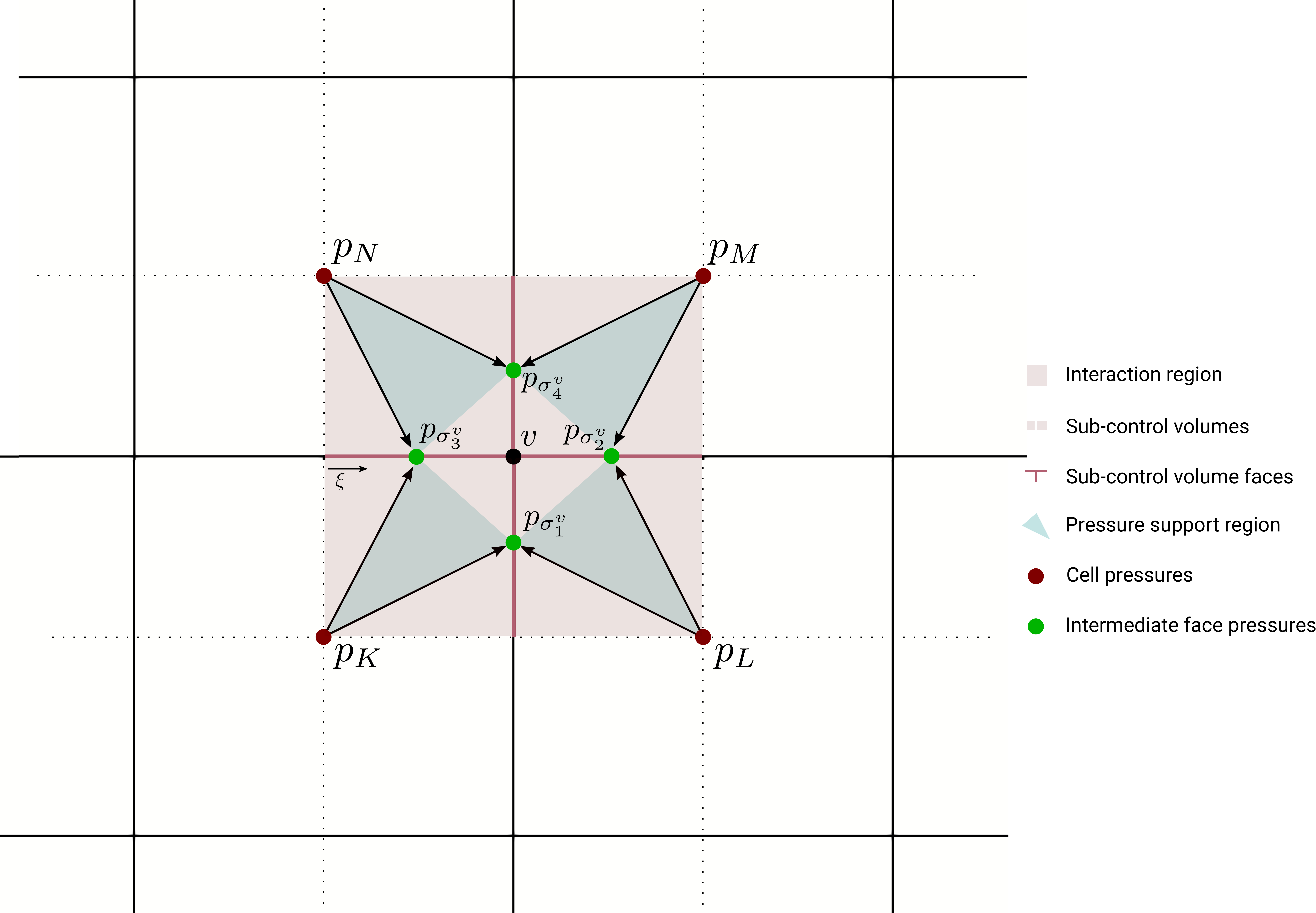}
  \caption{Interaction region for the MPFA-O method. The parameter $\xi$, $0 \le \xi < 1$ is used to define the location of the intermediate 
            face pressure unknowns $p_{\sigma_i^\vertex}$. Here, the situation for $\xi = 0.5$ is illustrated.}
  \label{pc:interactionRegion_interior}
\end{figure}

For the computation of the fluxes, a dual grid is created
 by connecting the barycenters of the cells with the barycenters of the faces ($d=2$)
 or the barycenters of the faces and edges ($d=3$). This divides each cell into sub-control
 volumes $K^\vertex$ with $K \in \mathcal{T}$ and $\vertex \in \calV$. Analogously, each face is sub-divided into sub-control volume faces $\sigma^\vertex$,
 see Figure \ref{pc:interactionRegion_interior}. Expressions for
the face fluxes $F_{K, \sigma^\vertex}$ are obtained by introducing the face
pressures $p_{\sigma^\vertex}$. The location of these face pressures along the sub-control volume faces $\sigma^\vertex$ 
is parameterized by $\xi$, $0 \le \xi < 1$, and is the center of the original face $\sigma$ for 
$\xi = 0$ and would be the position of the vertex $v$ for $\xi = 1$. These face pressures are then eliminated by enforcing
the continuity of fluxes across each sub-control volume face.
I.e., for each face $\sigma^\vertex$ between $K^\vertex$ and $L^\vertex$, we  impose:
\begin{subequations}
 \begin{align}
     &F_{K, \sigma^\vertex} + F_{L, \sigma^\vertex} = 0.
     \label{eq:sigmaConditions}
 \end{align}
\end{subequations}

We allow for piecewise constant $\mathbf{K}$
 (denoted as $\mathbf{K}_K$ for each cell $K$) and construct discrete gradients
 $\nabla_\mathcal{D}^{K^\vertex} p$, per sub-control volume $K^\vertex$, depending on its two embedded
 sub-control volume faces.
 Let us consider $K^\vertex$ in Figure~\ref{pc:interactionRegion_interior} with faces $\sigma_1^\vertex$ and $\sigma_3^\vertex$. Here, the discrete
 gradients are constructed to be consistent such that the following holds for $i \in \{1, 3\}$:
\begin{equation}
  \nabla_\mathcal{D}^{K^\vertex} p \cdot (\mathbf{x}_{\sigma^\vertex_i}- \mathbf{x}_{K}) =
  p_{\sigma^\vertex_i} - p_K.
\label{eq:piecewiselin}
\end{equation}

Thus, a discrete gradient (for sub-control volume $K^\vertex$) that fulfills the
two conditions \eqref{eq:piecewiselin} is defined by
\begin{equation}
  \nabla_\mathcal{D}^{K^\vertex} p  = \mathbb{D}^{-T}_{K^\vertex}
    \begin{bmatrix}
      p_{\sigma^\vertex_1} - p_K \\
      p_{\sigma^\vertex_3} - p_K
    \end{bmatrix}, \qquad \text{ with }\; \mathbb{D}_{K^\vertex} \eqbydef
    \begin{bmatrix}
      \mathbf{x}_{\sigma^\vertex_1}- \mathbf{x}_K & \mathbf{x}_{\sigma^\vertex_3} - \mathbf{x}_K
    \end{bmatrix}.
    \label{eq:MPFAGradientRecons}
\end{equation}

This enables us to write the discrete flux across $\sigma^\vertex$ between $K^\vertex$ and $L^\vertex$ as follows:
\begin{equation}
    F_{K, \sigma^\vertex} \eqbydef - |\sigma^\vertex| \mathbf{n}_{K, \sigma^\vertex}^T \mathbf{\mathbf{K}}_K \nabla_\mathcal{D}^{K^\vertex} p + \gamma_{K, \sigma^\vertex},
    \label{eq:discreteFlux}
\end{equation}
where we introduced $\gamma_{K, \sigma^\vertex} = \rho_{\sigma^\vertex} \meas{\sigma^\vertex} \mathbf{n}_{K, \sigma^\vertex}^T \mathbf{\mathbf{K}}_K \mathbf{g}$,
with $\rho_{\sigma^\vertex} = \frac{\rho_K+\rho_L}{2}$, to incorporate the effect of gravity.

To deduce a cell-centered scheme, the face pressures $p_{\sigma^\vertex_i}$
 are eliminated. This is done by enforcing flux continuity \eqref{eq:sigmaConditions}
 within each interaction
 volume and by solving a local system of equations. 

We rewrite these conditions in matrix form, and introduce the sans serif font to denote the corresponding matrices and vectors. All local
 face pressures of an interaction region are collected in the vector $\mathsf{p_{\sigma}}$, cell pressures in the vector
 $\mathsf{p_K}$, and all terms related to gravity in the vector $\mathsf{g}$. Flux continuity then allows us to rewrite the face pressures in terms of the cell pressures:
\begin{equation}
  \mathsf{A} \, \mathsf{p_{\sigma}} = \mathsf{B} \, \mathsf{p_K} + \Delta \mathsf{g}.
  \label{eq:ivLocalSystem}
\end{equation}

Here, the $\Delta$ represents the difference in contributions due to gravity over each face.
Let $\mathsf{f}$ denote the vector of all fluxes across the sub-control volume faces of the interaction
 region. These can be expressed in matrix form using equation \eqref{eq:discreteFlux}:
\begin{equation}
  \mathsf{f} = \mathsf{C} \, \mathsf{p_{\sigma}} + \mathsf{D} \, \mathsf{p_K} + \mathsf{g}.
  \label{eq:ivLocalFluxes}
\end{equation}

With these introduced matrices the final expressions for the local
 sub-control volume face fluxes, related to the interaction region, read:
\begin{equation}
  \mathsf{f} = \underbrace{\left( \mathsf{C} \mathsf{A}^{-1} \mathsf{B}
         + \mathsf{D} \right)}_{=: \mathsf{T}} \, \mathsf{p_K}
         + \mathsf{C} \mathsf{A}^{-1} \left( \Delta \mathsf{g} \right) + \mathsf{g}.
  \label{eq:ivLocalFluxesTrans}
\end{equation}
The entries of the matrix $\mathsf{T}$ are often referred to as the transmissibilities.


\subsection{Coupling}
\label{sub:discrete_coupling}

\begin{figure} [ht!]
  \centering
  \includegraphics[width=1.0\linewidth,keepaspectratio]{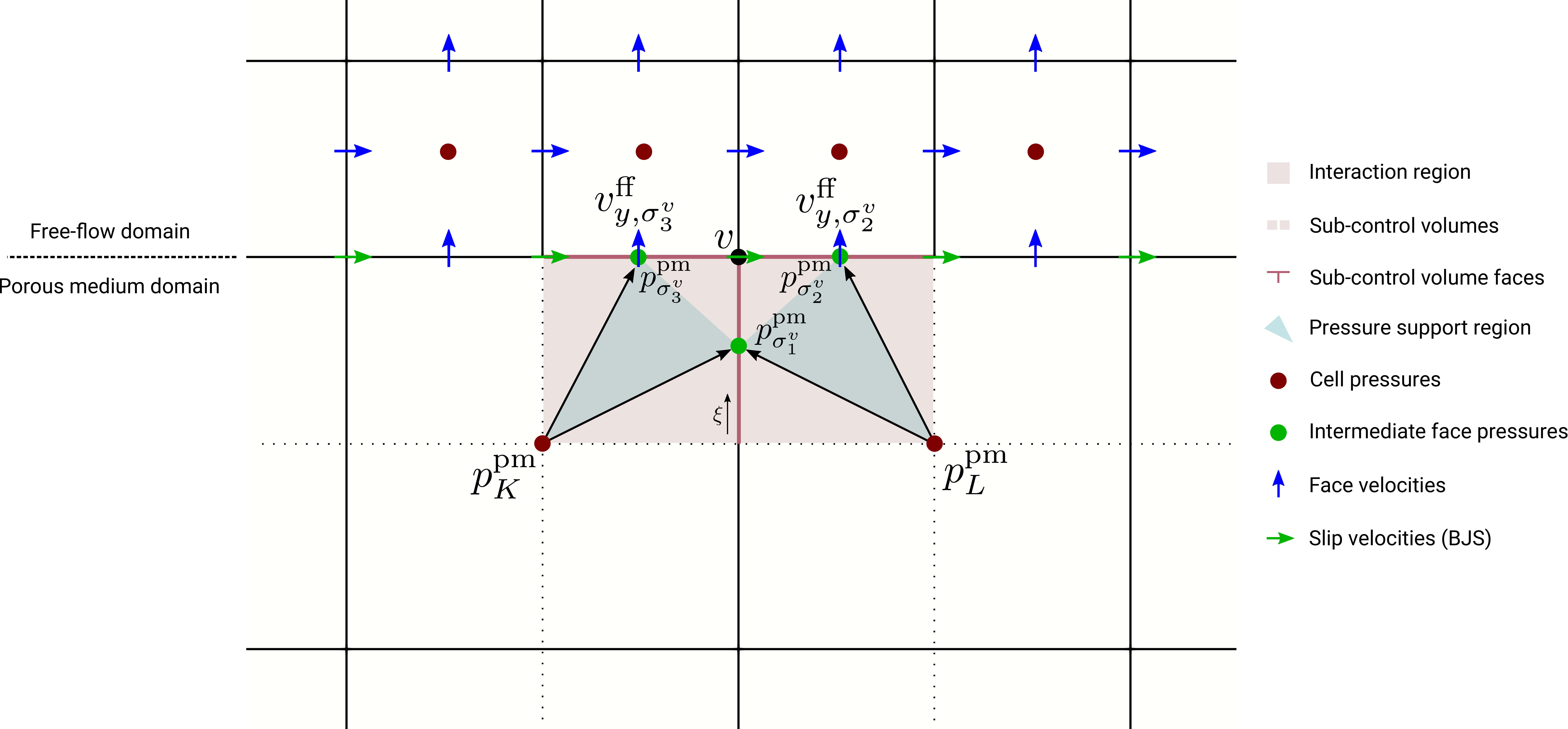}
  \caption{Interaction region for the MPFA-O method at the interface to the free-flow subdomain. The graphic illustrates the 
            non-matching grids at the interface and the choice of $\xi=0.5$ for the MPFA scheme such that the degrees of freedom 
            for the face velocities in the free-flow domain coincide with the intermediate face pressure unknowns introduced on the interface.}
  \label{pc:interactionRegion_interface}
\end{figure}

In this section, we consider the realization of the coupling conditions
\eqref{eq:interfaceConditions}. As depicted in Figure~\ref{pc:interactionRegion_interface}, the
grids are chosen to be non-matching at the interface such that each sub-control volume face
coincides with a face $\sigma \in \calE^\ff$ of the free-flow domain. In turn, a natural coupling
arises between the staggered grid discretization and the MPFA method, due to the coinciding degrees
of freedom (for $\xi = 0.5$, see Figure \ref{pc:interactionRegion_interface}).
 We emphasize that each cell in the porous-medium subdomain has two neighboring cells in
the free-flow subdomain (two-dimensional setup).

We start with the flux continuity condition \eqref{eq:interfaceConditionsFlux}. As depicted in
Figure~\ref{pc:interactionRegion_interface}, let us consider a sub-control volume $K^\vertex$ located at the
interface such that $\sigma_3^\vertex$ is located on $\Gamma^\ipmff$. We then let the velocity from the free-flow
 domain determine the flux over $\sigma_3^\vertex$:
\begin{equation}
\begin{aligned}
    F_{K, \sigma_3^\vertex} &= - \mu^{up} \meas{\sigma_3^\vertex} \veli{y, \sigma_3^\vertex}^\ff.
  \label{eq:interactionRegionConditions}
\end{aligned}
\end{equation}
 
Collecting the right-hand side of \eqref{eq:interactionRegionConditions} in the matrix-vector product
 $\mathsf{Wv_\ff}$, we can rewrite \eqref{eq:ivLocalSystem} for interaction regions
 that are located at the interface to the free-flow domain:
\begin{equation}
  \mathsf{A} \, \mathsf{p_{\sigma}} = \mathsf{B} \, \mathsf{p_K} + \mathsf{Wv_\ff} + \Delta \mathsf{g}.
  \label{eq:ivLocalSystemInterface}
\end{equation}

In the situation shown in Figure~\ref{pc:interactionRegion_interface}, the face pressures
$p^{\porm}_{\sigma^\vertex_1}, p^{\porm}_{\sigma^\vertex_2}, p^{\porm}_{\sigma^\vertex_3}$ are thus dependent on the
primary unknowns $p^\porm_K, p^\porm_L$ of the porous-medium domain and the face velocities
$v^\ff_{\sigma^\vertex_2},v^\ff_{\sigma^\vertex_3}$ of the free-flow domain (i.e. $p^{\porm}_{\sigma^\vertex_i} =
p^{\porm}_{\sigma^\vertex_i}(p^\porm_K, p^\porm_L,v^\ff_{\sigma^\vertex_2},v^\ff_{\sigma^\vertex_3})$). Insertion of
\eqref{eq:ivLocalSystemInterface} in \eqref{eq:ivLocalFluxes} leads to the expression
\begin{equation}
  \mathsf{f} = \mathsf{T} \, \mathsf{p_K}
         + \mathsf{C} \mathsf{A}^{-1} \left( \mathsf{Wv_\ff} + \Delta \mathsf{g} \right) + \mathsf{g}
  \label{eq:ivLocalFluxesInterface}
\end{equation}
for the fluxes across sub-control volume faces within interaction regions located at the interface to the
 free-flow domain.

The remaining coupling conditions are imposed as followed. The momentum balance
\eqref{eq:interfaceConditionsMomentum} is enforced using the reconstructed face pressures from
\eqref{eq:ivLocalSystemInterface}. On the other hand, the Beavers-Joseph-Saffman condition
\eqref{eq:interfaceConditionsBJS} is technically a boundary condition for the free-flow problem, as
previously noted in Section~\ref{sec:equations}, and is implemented accordingly.

Finally, we remark that in the case of compressible fluids, the density and viscosity are
pressure-dependent. The upwind terms $\mu^{\porm,\up}, \varrho^{\porm,\up}, \varrho^{\ff,\up}$ are
then evaluated using the cell-pressure unknowns, i.e. for a face $\sigma$ between porous-medium
cell $K$ and free-flow cell $L$, these terms therefore depend on the pressures $p^\porm_K$ and
$p^\ff_L$.

\section{Numerical results}
\label{sec:numresults}
All simulations are performed using the open-source simulator \Dumux \cite{dumux}, which comes in the form of an additional DUNE module \citep{blatt.ea:2016}.
We employ a monolithic approach, where both sub-problems are assembled into one system of equations and use an implicit Euler method for the time discretization.
Newton's method is used to solve the non-linearities involved in the systems of equations.
For all test cases, the compressible fluid ``air'' (see the \Dumux documentation \cite{dumux}) is used. We consider two-dimensional setups, however, the implementation is also able
to handle three-dimensional domains.

\subsection{Test case 1}

The first test case is similar to the one that has been presented in \cite{iliev2004a}. 
However, the authors of \cite{iliev2004a} used a Navier-Stokes-Brinkman-type system for 
both domains and this system is discretized by using a staggered-grid scheme in both domains. 
For anisotropic permeability tensors, this requires the interpolation of all velocity components at grid faces. 

In this work, we use a different approach, where a staggered-grid scheme is used to discretize 
the free-flow system and a cell-centered finite volume scheme to discretize the porous-medium system.
Thus, no additional velocity degrees of freedom are needed for the porous-medium domain. 
However, for anisotropic permeability tensors or unstructured grids this requires more sophisticated cell-centered finite volume schemes, 
as for example the MPFA scheme that has been presented in Section \ref{sub:ccfvm}.

\begin{figure} [ht!]
	\centering
	\includegraphics[width=0.7\linewidth,keepaspectratio]{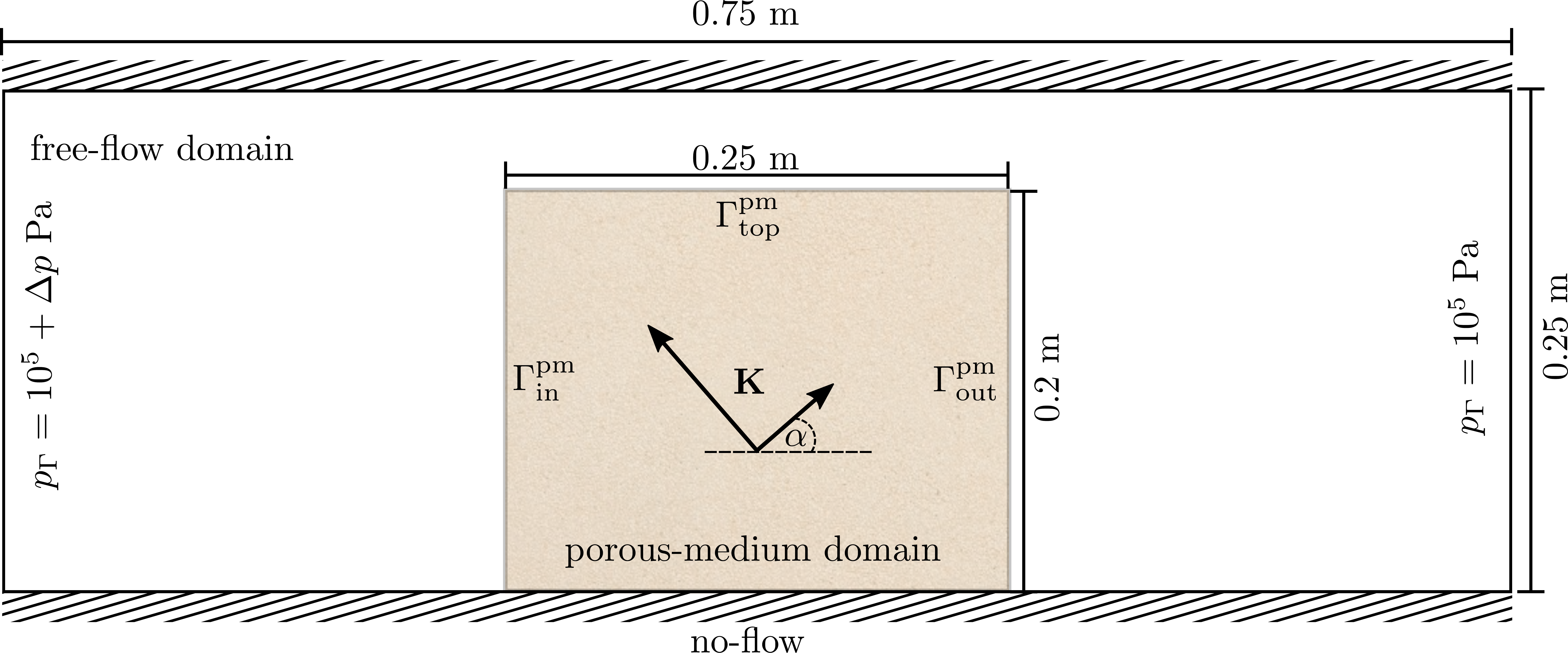}
	\caption{Setting for the first test case.}
	\label{pc:settingCaseOne}
\end{figure}

Air is flowing through a two-dimensional channel which is partially blocked by a rectangular porous medium as shown in Figure \ref{pc:settingCaseOne}. The first test case involves small Reynolds
numbers ($Re \ll 1$) with respect to the average velocity in the narrow section in the channel above the porous medium. In this case, only the stationary solution (which is reached after a few time
steps, starting from a resting fluid) is investigated. The top and bottom of the domain are considered as rigid, impermeable walls with $\vel = 0$ (including the wall part below the porous box). Flow is driven by a pressure difference between the left and the right
boundary which is set to $\Delta p = \unit[10^{-6}]{Pa}$.

The permeability tensor in the porous medium is given as
\begin{equation}
\mathbf{K} =  \mathbf{R}(\alpha)
\begin{pmatrix}
\frac{1}{\beta}k & 0 \\
0 & k \\
\end{pmatrix} \mathbf{R}^{-1}(\alpha), \quad
\text{with} \quad \mathbf{R}(\alpha) = \begin{pmatrix}
\cos \alpha & -\sin \alpha \\
\sin \alpha & \cos \alpha \\
\end{pmatrix},
\label{eq:RotatedPermTensor}
\end{equation}
where $\alpha$ is the rotation angle, $\beta \geq 1$ the anisotropy ratio, and $\unit[k=10^{-6}]{m^2}$. In the following, the influence of the anisotropy on the total mass fluxes crossing the boundaries $\Gamma^\porm_\mathrm{in}, \Gamma^\porm_\mathrm{out}, \Gamma^\porm_\mathrm{top}$ is investigated for $\alpha \in \lbrace  -45^\circ, -30^\circ, 0^\circ, 30^\circ, 45^\circ  \rbrace$ and $\beta \in \lbrace  10 ,100  \rbrace$. 
In the free-flow domain the same grid is used for the MPFA and the TPFA schemes, whereas in the porous medium the grid used for the MPFA scheme is coarsened such that we obtain a non-matching interface, as shown in Figure \ref{pc:interactionRegion_interface}.

\begin{figure}[ht!]
	\centering
	\includegraphics[width=0.45\linewidth,keepaspectratio]{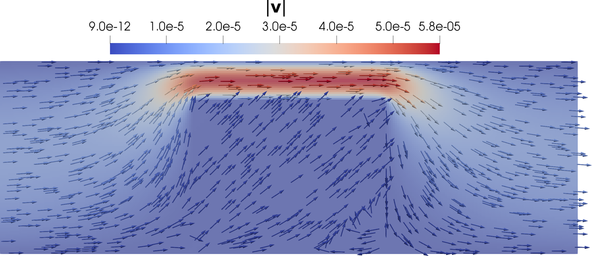}
	\hspace{0.04\linewidth}
	\includegraphics[width=0.45\linewidth,keepaspectratio]{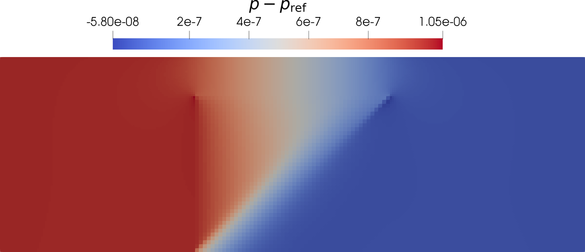} \\
	\vspace{0.02\linewidth}
	\includegraphics[width=0.45\linewidth,keepaspectratio]{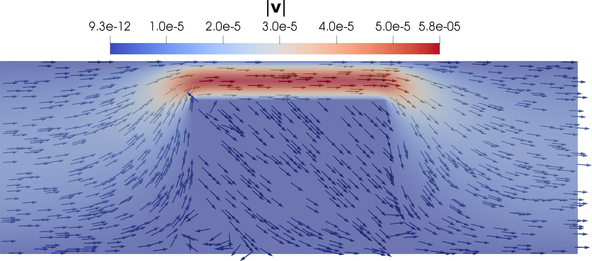}
	\hspace{0.04\linewidth}
	\includegraphics[width=0.45\linewidth,keepaspectratio]{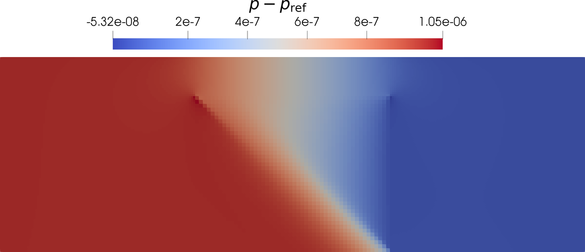}
	\caption{Velocity (left column) and pressure (right column) profiles of MPFA scheme for an ansisotropy ratio of $\beta = 100$ and for angles $\alpha = -45^\circ$ (upper row) and $\alpha = 45^\circ$ (lower row). The reference pressure is set to $p_{\mathrm{ref}} = \unit[10^5]{Pa}$.}
	\label{pc:testOneVelMpfa}
\end{figure}

\begin{figure}[ht!]
	\centering
	\includegraphics[width=0.45\linewidth,keepaspectratio]{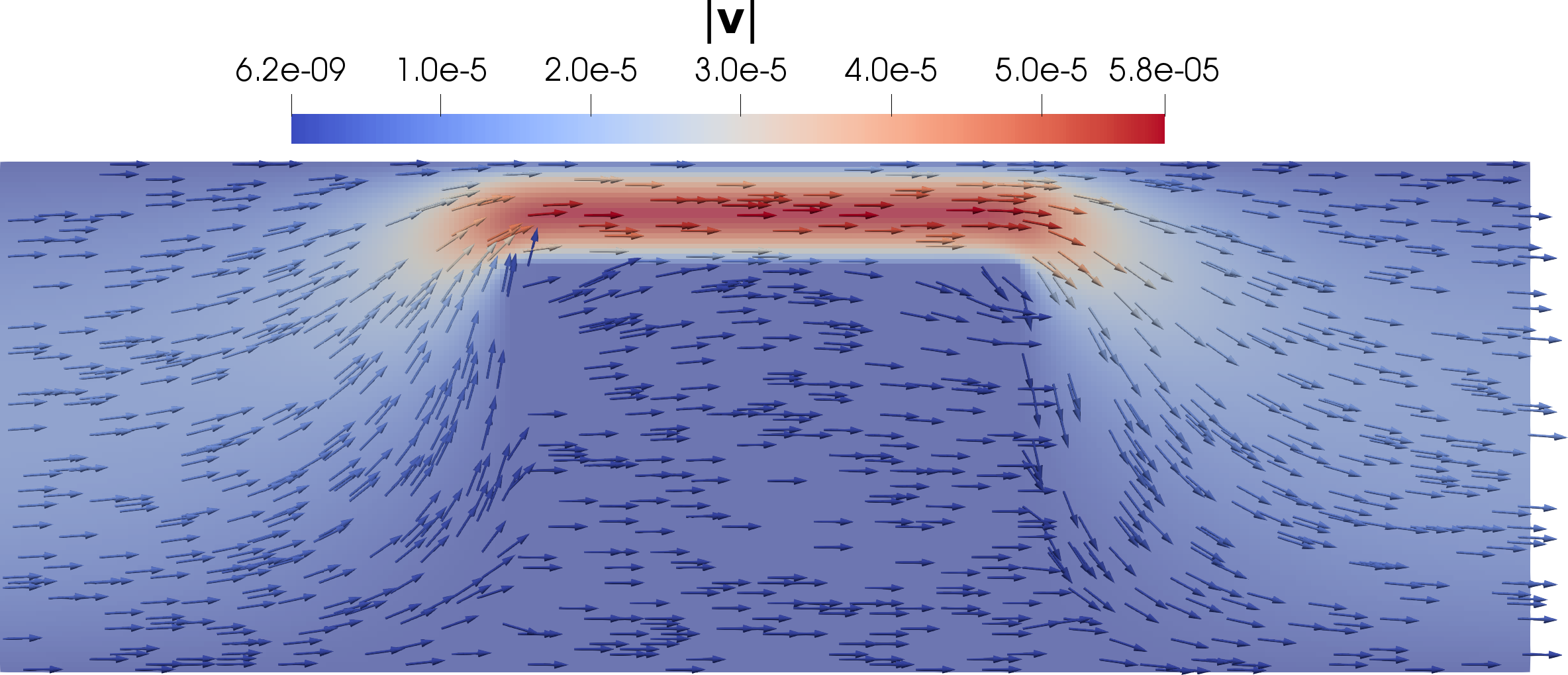}
	\hspace{0.04\linewidth}
	\includegraphics[width=0.45\linewidth,keepaspectratio]{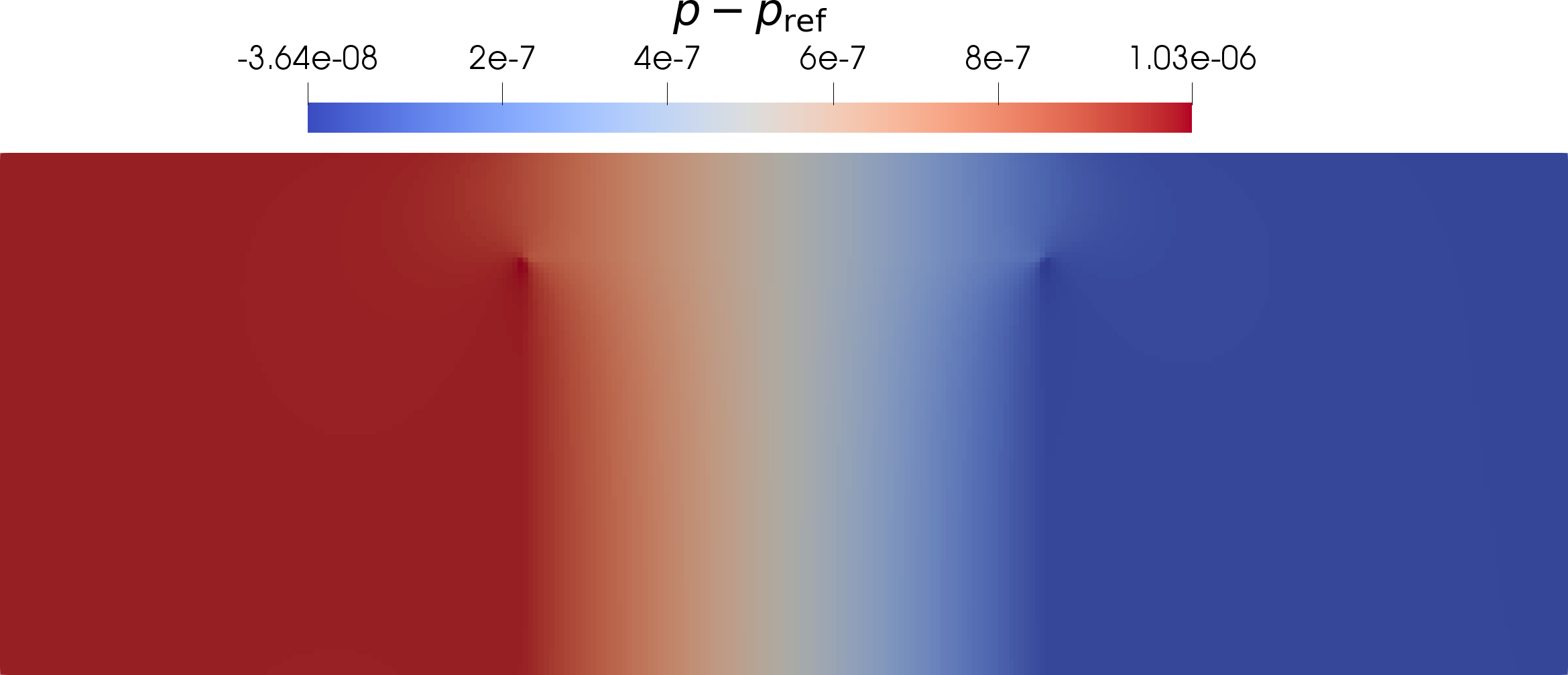}
	\caption{Velocity (left column) and pressure (right column) profiles of TPFA scheme for an ansisotropy ratio of $\beta = 100$ and for angle $\alpha = 45^\circ$. The reference pressure is set to $p_{\mathrm{ref}} = \unit[10^5]{Pa}$.}
	\label{pc:testOneVelTpfa}
\end{figure}

Figure \ref{pc:testOneVelMpfa} and \ref{pc:testOneVelTpfa} show the velocity magnitude and the pressure profile
of the MPFA and TPFA schemes for $\beta = 100$ and $\alpha \in \lbrace -45^\circ, 45^\circ\rbrace$.
For the TPFA scheme, the results are only shown for $\alpha = 45^\circ$ because the results for $\alpha = - 45^\circ$ are identical.
Due to the flow resistance imposed by the porous medium, most of the air passes this obstacle through the constricted section above the block,
thus leading to the highest flow velocities there. While the gas passes the porous block virtually parallel for both $\alpha = 45^\circ$ and $\alpha = 0^\circ$
when applying the TPFA method, the effect of anisotropy is clearly visible in the MPFA results. Here, the flow follows the inclined principal direction
of the permeability tensor, exiting or entering the porous domain at the top. Small regions of local recirculation can be found within the obstacle which are caused by the medium's
anisotropy and the closed wall at the bottom. For $\alpha = - 45^\circ$, the upward flow in the box creates a recirculation at the right part of the bottom, where a small amount of
gas is actually pulled from the free-flow channel into the porous medium. Analogously for $\alpha = 45^\circ$, the downward flow causes a recirculation at the bottom left, where the
gas cannot exit through the solid wall and thus leaves the domain towards the left, in opposition to the general flow field.

Due to the small pressure differences, we subtract a reference value
of $\unit[p_\mathrm{ref}=10^{5}]{Pa}$ for improved visibility as shown on the right side of Figures \ref{pc:testOneVelMpfa} and \ref{pc:testOneVelTpfa}.
Almost the entire pressure drop along the channel is observed at the porous domain.
Again, the influence of anisotropy is clearly visibly for the MPFA results in terms of an inclined pressure profile which also reflects the closed bottom of the domain. None of these effects can be
captured by the TPFA method.

\begin{figure}[ht!]
	\centering
	\includegraphics[width=0.45\linewidth,keepaspectratio]{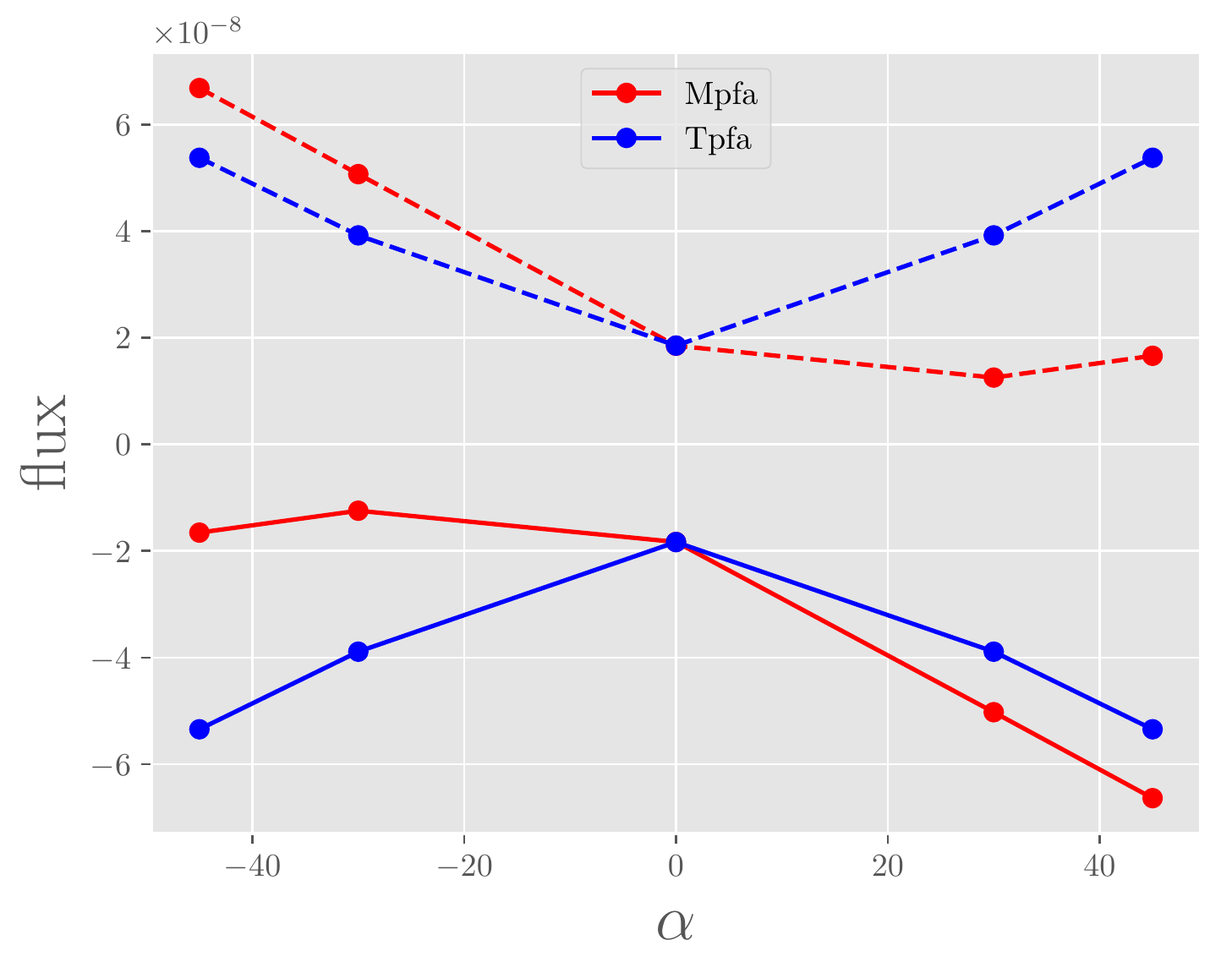}
	\hspace{0.04\linewidth}
	\includegraphics[width=0.45\linewidth,keepaspectratio]{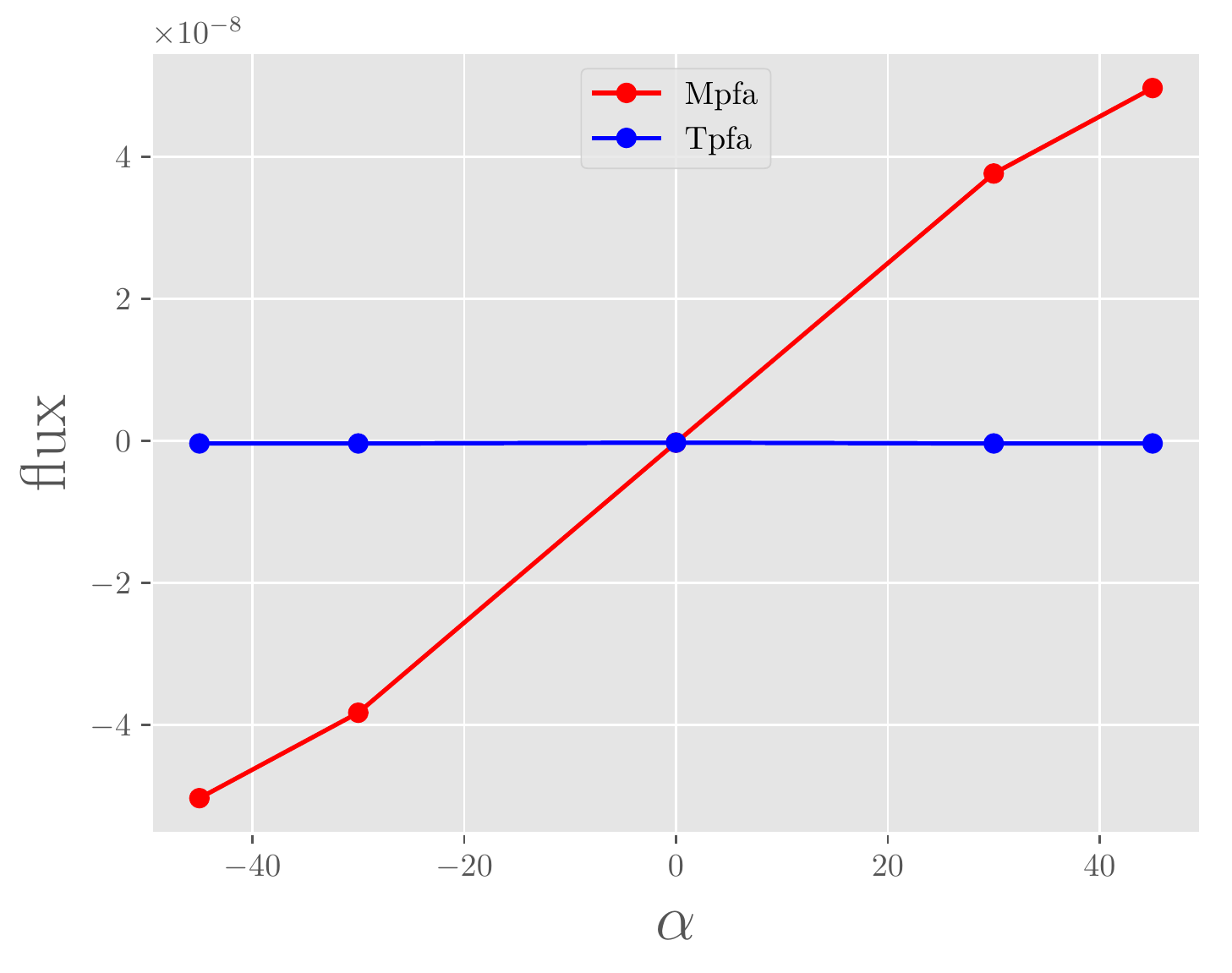}
    \includegraphics[width=0.45\linewidth,keepaspectratio]{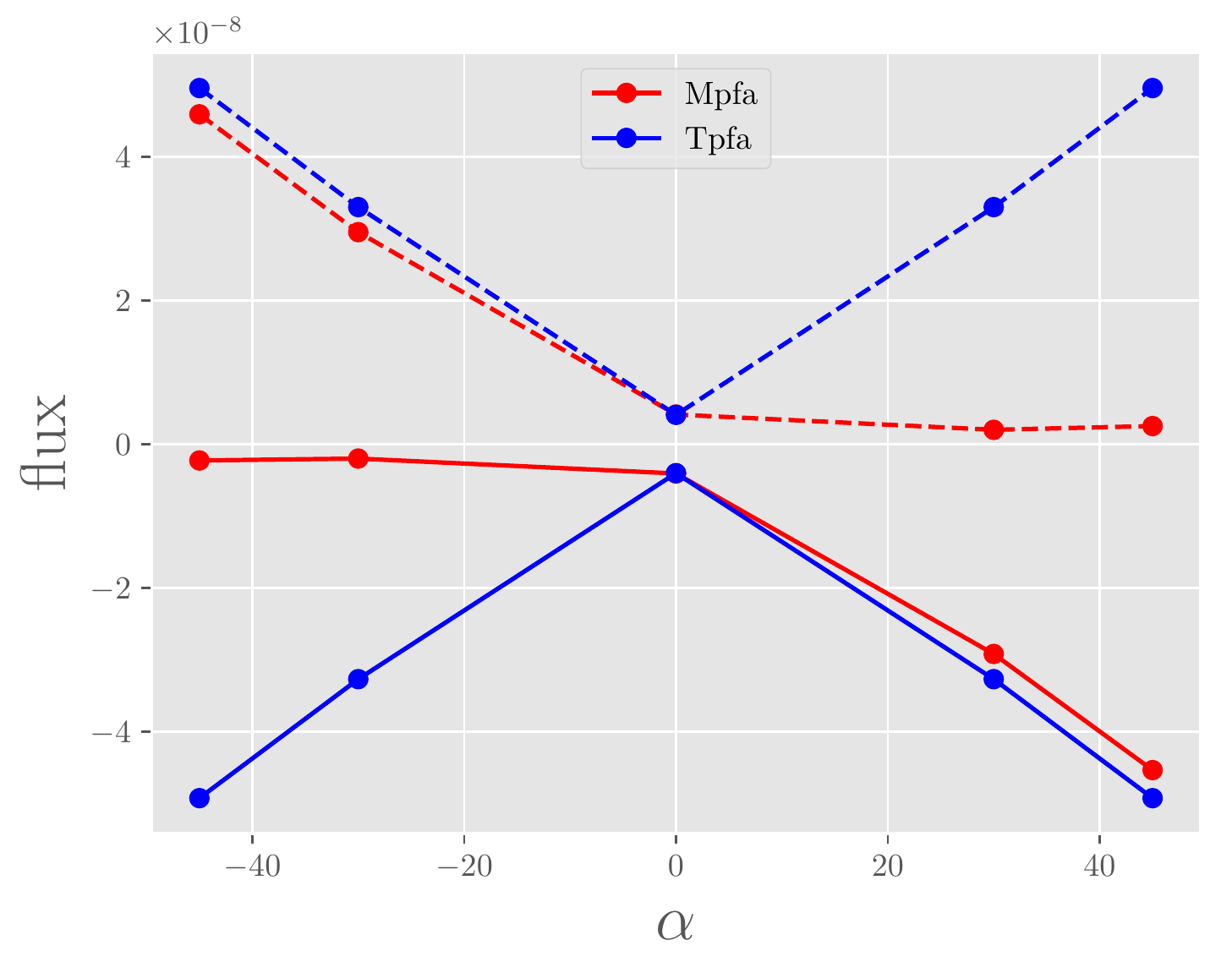}
	\hspace{0.04\linewidth}
	\includegraphics[width=0.45\linewidth,keepaspectratio]{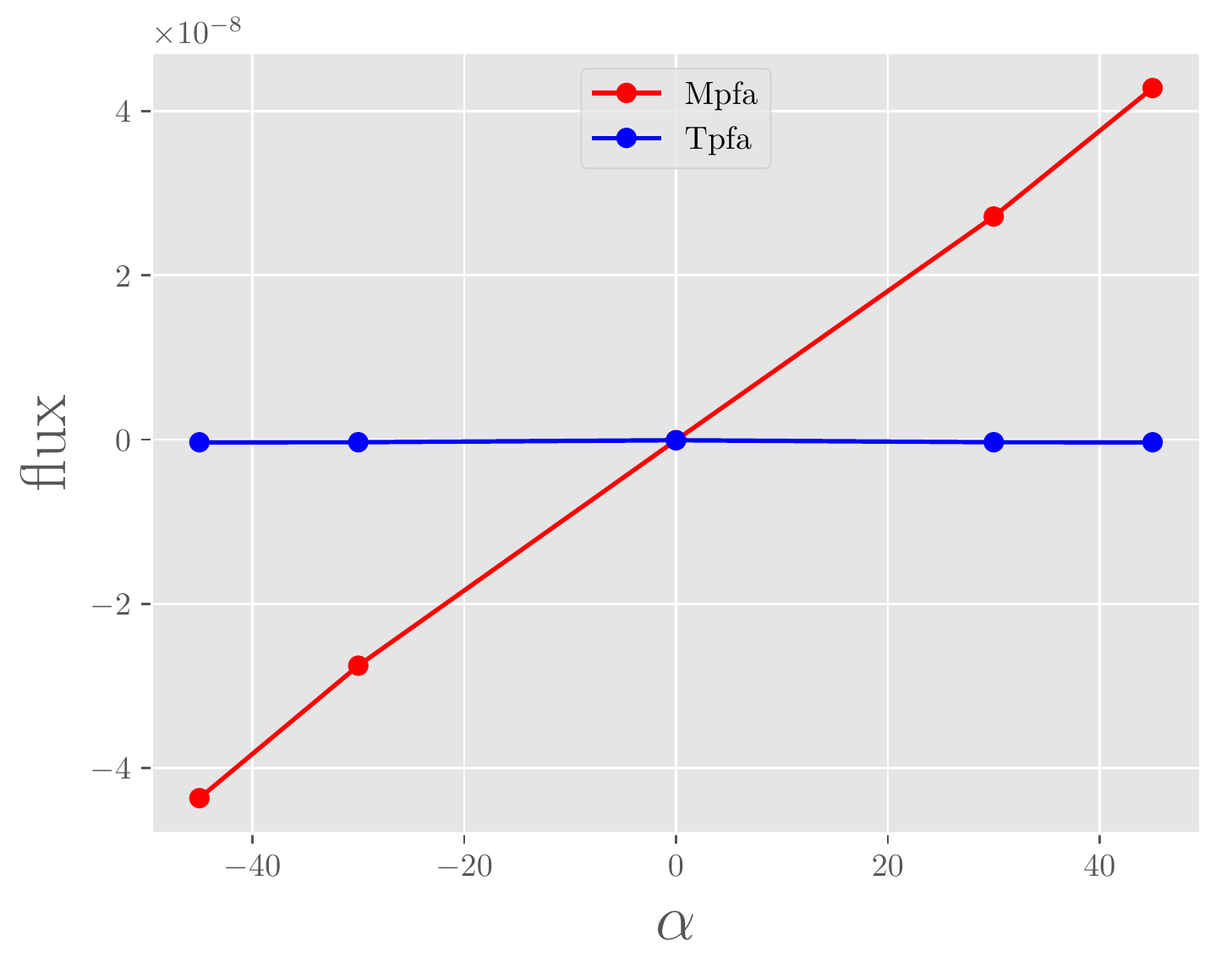}
	\caption{Total mass fluxes of MPFA and TPFA schemes over the porous-medium boundaries $\Gamma^\porm_\mathrm{in}, \Gamma^\porm_\mathrm{out}$ (left column) and  $\Gamma^\porm_\mathrm{top}$ (right column) for an ansisotropy ratio of $\beta = 10$ (upper row) and $\beta = 100$ (lower row). In the left pictures, the dashed lines correspond to the mass fluxes at the inlet boundary $\Gamma^\porm_\mathrm{in}$, whereas the solid lines to the fluxes at the outlet boundary $\Gamma^\porm_\mathrm{out}$. Positive fluxes mean that fluid flows into the porous medium.}
	\label{pc:massFluxesTestOne}
\end{figure}

Figure \ref{pc:massFluxesTestOne} shows the total mass fluxes over the porous-medium free-flow interfaces. For $\alpha = 0$, the TPFA and MPFA scheme result in almost
the same solutions (small differences occur because a finer mesh is used for the TPFA scheme in the porous medium). For this angle, the TPFA scheme is also consistent
and produces the correct results. However, it is observed that for all other angles the total mass fluxes significantly differ. The off-diagonal terms of $\mathbf{K}$
are not considered in the TPFA transmissibilities, which explains why the TPFA results are independent from the direction of rotation. Furthermore, the total fluxes
over $\Gamma^\porm_\mathrm{top}$ are small for the TPFA scheme, in contrast to the MPFA scheme, where the total mass fluxes increase with increasing rotation angle
due to the contribution from the non-parallel porous-medium flow as mentioned above.

\subsection{Test case 2}
The next test case investigates the solution behavior for higher Reynolds numbers.
It uses a similar setting as the previous one with the difference that the channel is elongated in $x$-direction. The computational domains are given by $\Omega = [0,2.5] \times [0,0.25]$ m and $\Omega^\porm = [0.4,0.6] \times [0,0.2]$ m such that $\Omega^\ff = \Omega \setminus \Omega^\porm$. A pressure difference of $\Delta p = \unit[2 \cdot 10^{-3}]{Pa}$ between the left and the right boundary results in $Re \approx 130$ in the channel right atop the porous block.

\begin{figure}[ht!]
	\centering
	\includegraphics[width=0.7\linewidth,keepaspectratio]{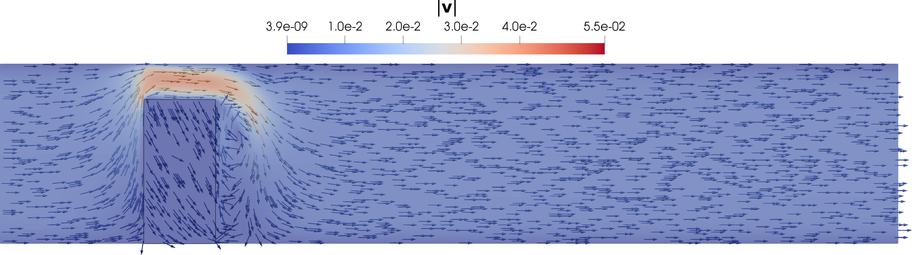}
	\includegraphics[width=0.7\linewidth,keepaspectratio]{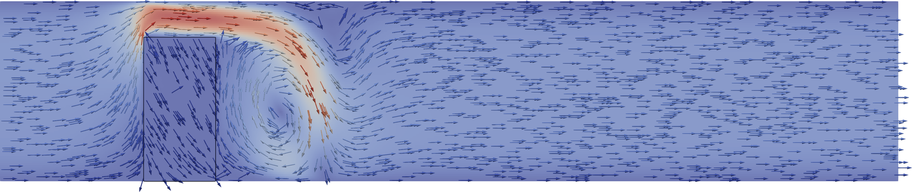}
	\includegraphics[width=0.7\linewidth,keepaspectratio]{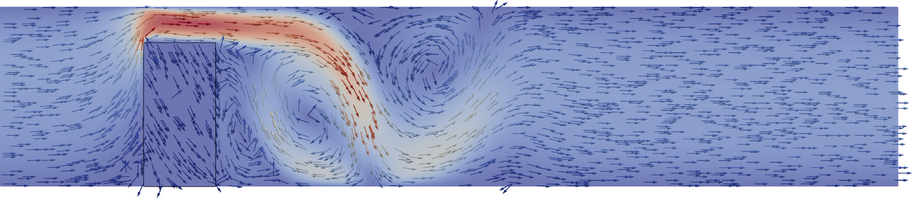}
	\includegraphics[width=0.7\linewidth,keepaspectratio]{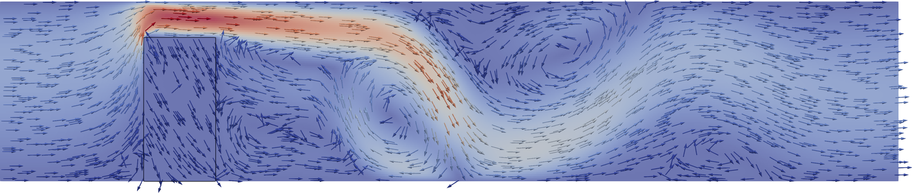}
	\includegraphics[width=0.7\linewidth,keepaspectratio]{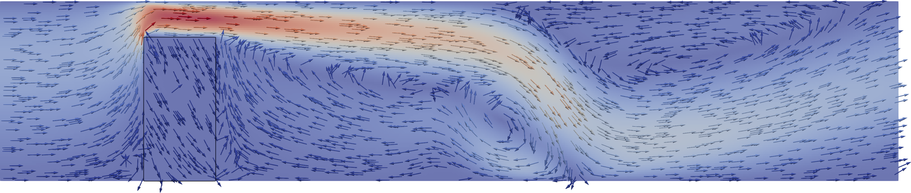}
	\caption{Velocity profiles of MPFA scheme for an ansisotropy ratio of $\beta = 100$ and for angle $\alpha = 45^\circ$ for times $t \in \lbrace \unit[20]{s}, \unit[40]{s}, \unit[80]{s}, \unit[200]{s}, \unit[1000]{s} \rbrace$. The domain is scaled by a factor of 2 in $y$-direction for a better visualization. The porous-medium boundary is represented by the black lines. }
	\label{pc:testHighReVel}
\end{figure}

\begin{figure}[ht!]
	\centering
	\includegraphics[width=0.7\linewidth,keepaspectratio]{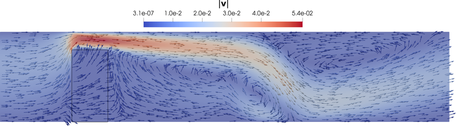}
	\caption{Velocity profile of TPFA scheme for an ansisotropy ratio of $\beta = 100$ and for angle $\alpha = 45^\circ$ for times $t = \unit[1000]{s}$. The domain is scaled by a factor of 2 in $y$-direction for a better visualization. The porous-medium boundary is represented by the black lines. }
	\label{pc:testHighReVelTPFA}
\end{figure}

\begin{figure}[ht!]
	\centering
	\includegraphics[width=0.45\linewidth,keepaspectratio]{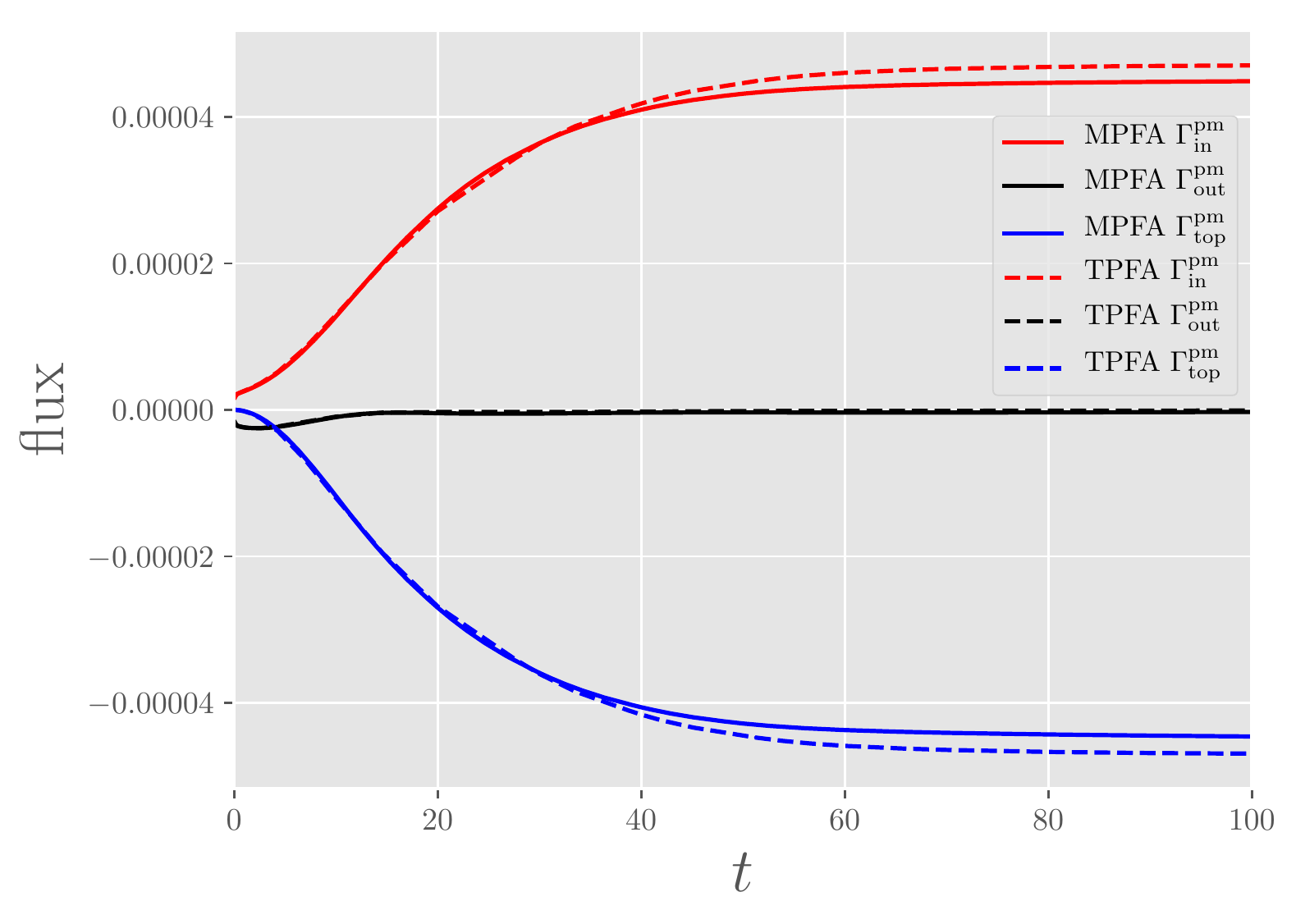}
	\hspace{0.04\linewidth}
    \includegraphics[width=0.45\linewidth,keepaspectratio]{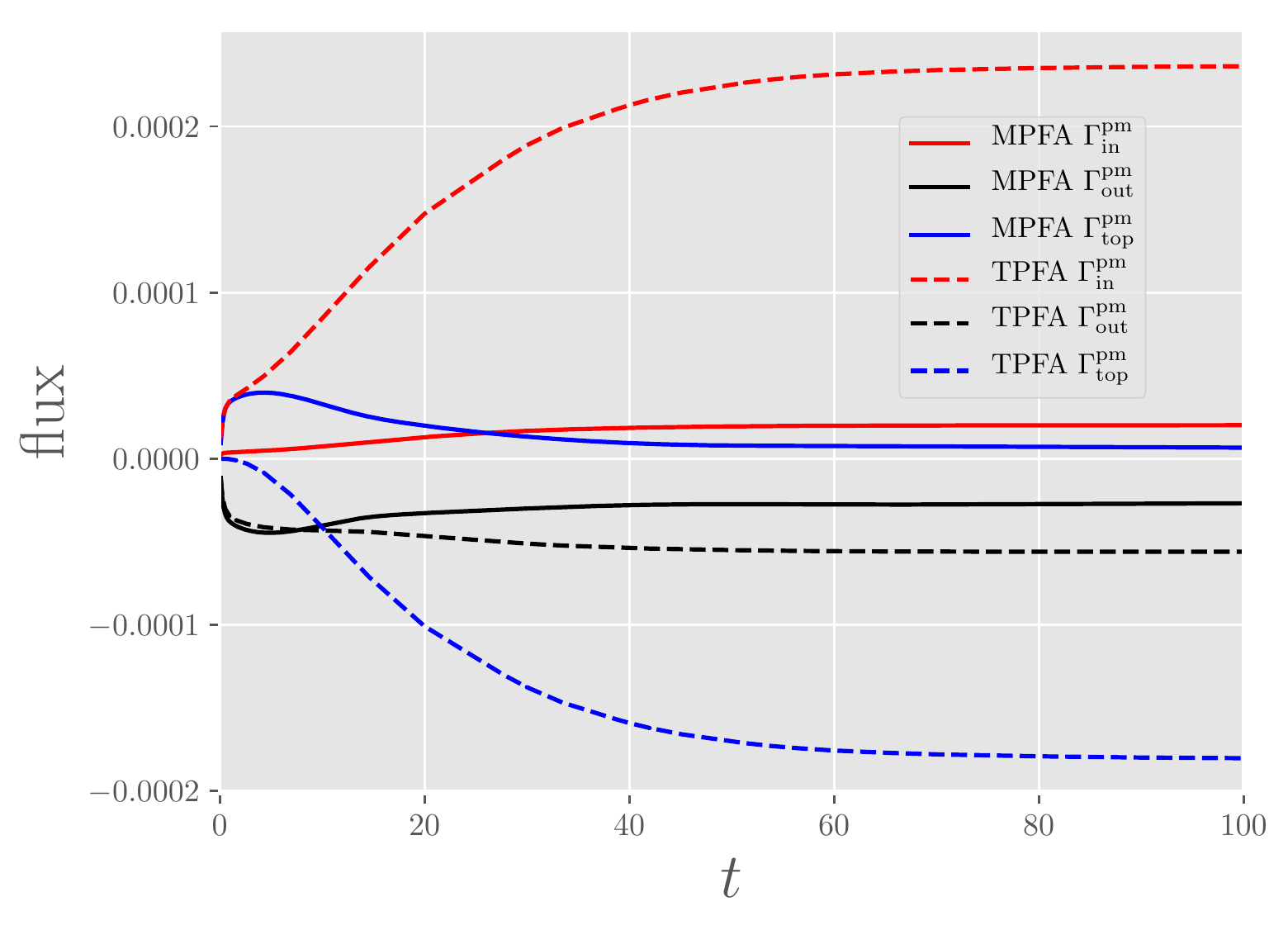}
	\includegraphics[width=0.44\linewidth,keepaspectratio]{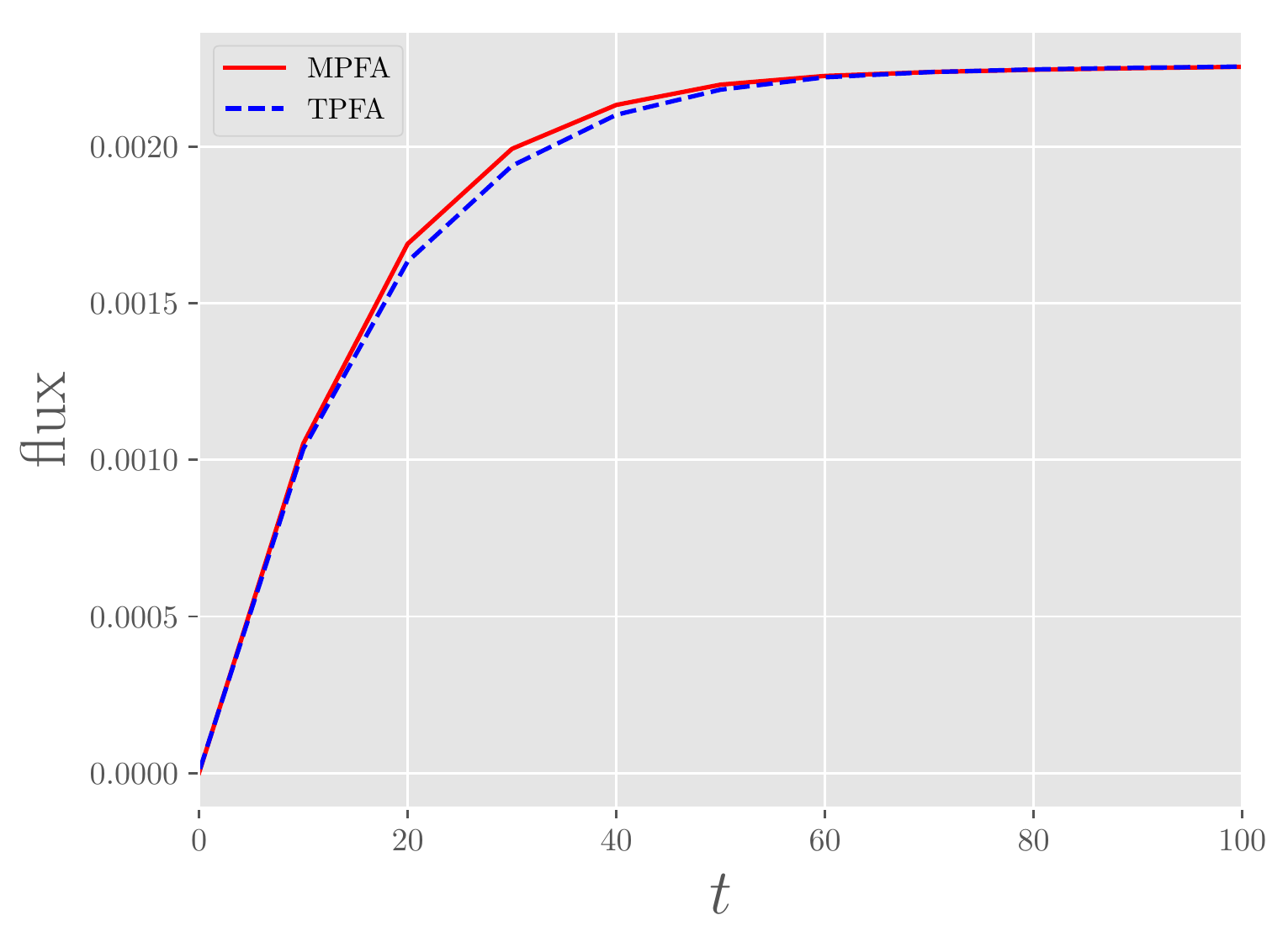}
	\hspace{0.06\linewidth}
	\includegraphics[width=0.44\linewidth,keepaspectratio]{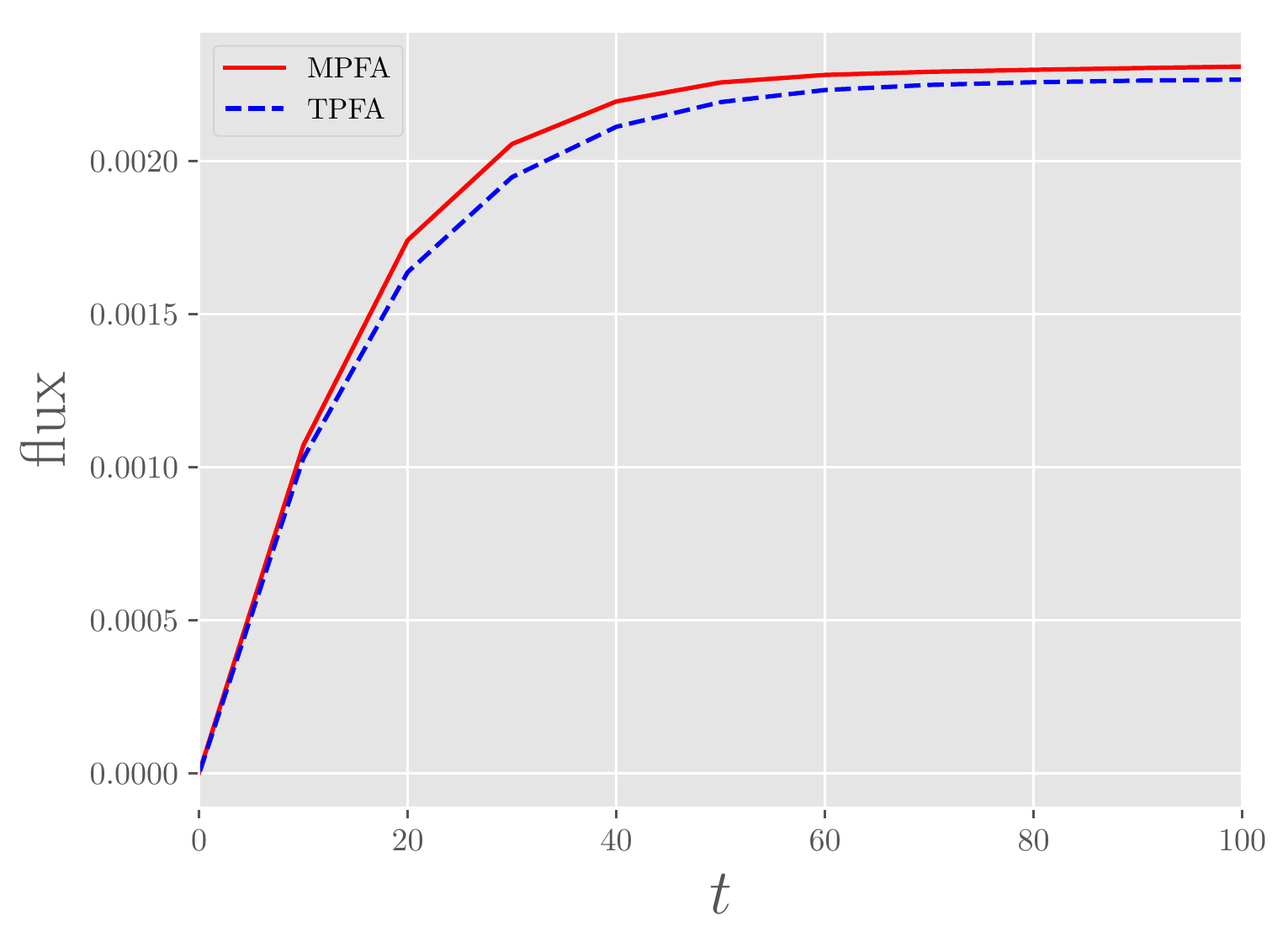}
	\caption{Total mass fluxes of MPFA and TPFA schemes over the porous-medium boundaries $\Gamma^\porm_\mathrm{in}, \Gamma^\porm_\mathrm{top}, \Gamma^\porm_\mathrm{out}$ (upper row) and  at the center of the constricted section of the free-flow channel (over the line segment connecting the points $(\unit[0.5]{m},\unit[0.2]{m})^T$ and $(\unit[0.5]{m},\unit[0.25]{m})^T$) (lower row) for an ansisotropy ratio of $\beta = 0$ (left column) and $\beta = 100$ (right column). Positive fluxes mean that fluid flows into the porous medium.}
	\label{pc:massFluxesTestTwo}
\end{figure}

For this test case, it takes much longer until a stationary solution is reached, which is why the solutions are investigated at different time steps.
In the following, we focus on the discussion of the case where $\alpha = 0^\circ$ or $\alpha = 45^\circ$ and $\beta = 100$. Figure \ref{pc:testHighReVel} shows the resulting velocity fields of the MPFA method at
different times for $\alpha = 45^\circ$. As before, the gas hits the front of the porous block and is forced mainly through the narrow channel section above it which acts as
some sort of duct through which a jet of high velocity fluid streams in the open channel on the right side. As the velocity gradually increases over time, the formation of vortex structures within
the free-flow channel can be observed after around $\unit[20]{s}$. Having reached equilibrium at around $\unit[1000]{s}$, the system features two stable countercurrent, larger vortices downstream
the porous obstacle and one small recirculation zone in front of the block. Again, the flow within the obstacle is clearly influenced by the anisotropy, a feature not reproduced by the TPFA method
as seen in Figure \ref{pc:testHighReVelTPFA}. Here, the flow within the block does not follow $\mathbf{K}$'s orientation. Instead, the fluid immediately strives for the top of the obstacle once it
has entered it which is due to the strongly increased vertical permeability ($\beta = 100$). Again, the TPFA scheme is not able to capture this relevant anisotropy effect as it only considers the
main diagonals of $\mathbf{K}$. Figure \ref{pc:massFluxesTestTwo} depicts the cumulative mass fluxes per unit depth across the porous medium's boundaries (top row) and across the plane at the center
of the narrow channel section ($x = \unit[0.5]{m}$, bottom row) for a case with $\alpha = 0^\circ$ (left column) and one with $\alpha = 45^\circ$ (right column). For the former, the results of
the TPFA and MPFA are very similar. With increasing time, the total mass flux entering the obstacle's front (red line) increases in accordance with the global velocity field (see Figure
\ref{pc:testHighReVel}) until it approaches a constant value after approximately $\unit[100]{s}$. The same applies for the fluxes leaving the obstacle's top (blue line), which approaches the same
value as the incoming fluxes but with a different sign. This indicates that no mass flux occurs over the obstacle's back which can also be seen from the black line. This line (and thus, the fluxes
over the obstacle's top) only deviates from zero at the very beginning of the simulation where it has to balance out the initial disequilibrium between the red and blue curve. Due to the imposed
anisotropy ratio of $\beta = 100$, the  horizontal permeability is 100 times lower than the vertical one and fluid is immediately pushed towards the top, where it exits the domain again.
Even though the fluid's density is actually pressure-dependent, significant compressibility effects could not be
observed for the given setup. The differences between TPFA and MPFA  with respect to the blue and red curves (front and back of the obstacle) are most likely due to the different
discretization width in the porous domain as explained before. In total, less gas seems to enter the porous block for the MPFA method. The plot on the lower left of Figure
\ref{pc:massFluxesTestTwo} shows the temporal evolution of the total mass flux through the channel above the obstacle. Both the TPFA and the MPFA scheme converge to the same result,
the grid resolution of the free-flow domain is identical for both cases. Note that even with a constant flux through the constriction after $\unit[100]{s}$, the global velocity field downstream the
obstacle still changes, including the formation of vortices as described before.

Considering the right column of Figure \ref{pc:massFluxesTestTwo} ($\alpha = 45^\circ$), the differences between the two schemes are significantly higher. For the MPFA method, we observe a
considerable inflow across the obstacle's top boundary, coming from the constricted free-flow channel and drawn by the downwards inclined flow field within the porous medium. At the same time,
there is only a very limited inflow through the obstacle's front which is due the obstacle's anistropy.

This is in strong contrast to the TPFA method's result, where air is still leaving the obstacle's top as also described before and seen in Figure \ref{pc:testHighReVelTPFA}.
The MPFA method results in a higher flux through the constricted channel (see bottom right of Figure \ref{pc:massFluxesTestTwo}) as there is less fluid entering the porous obstacle as explained
above.

\subsection{Test case 3}

\begin{figure}[ht!]
	\centering
	\includegraphics[width=0.95\linewidth,keepaspectratio]{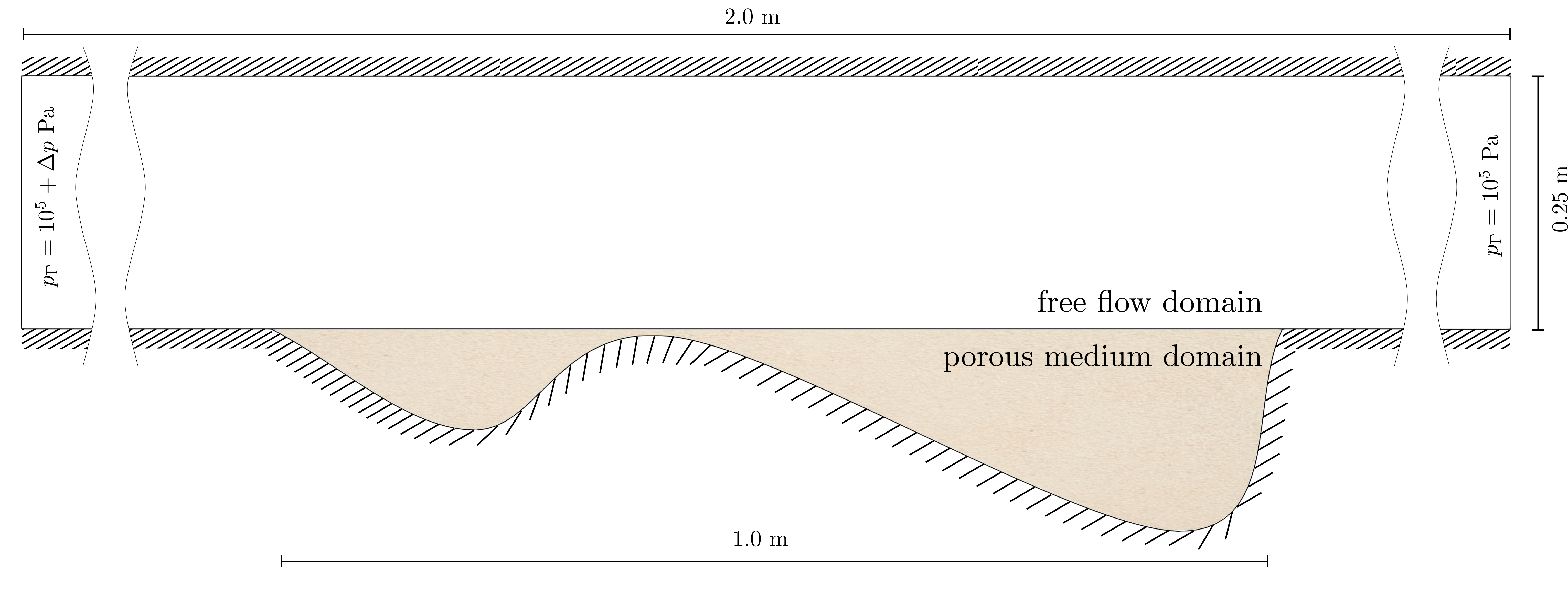}
	\caption{Setting for the third test case.}
	\label{pc:settingCaseTwo}
\end{figure}

\begin{figure}[ht!]
	\centering
	\includegraphics[width=0.55\linewidth,keepaspectratio]{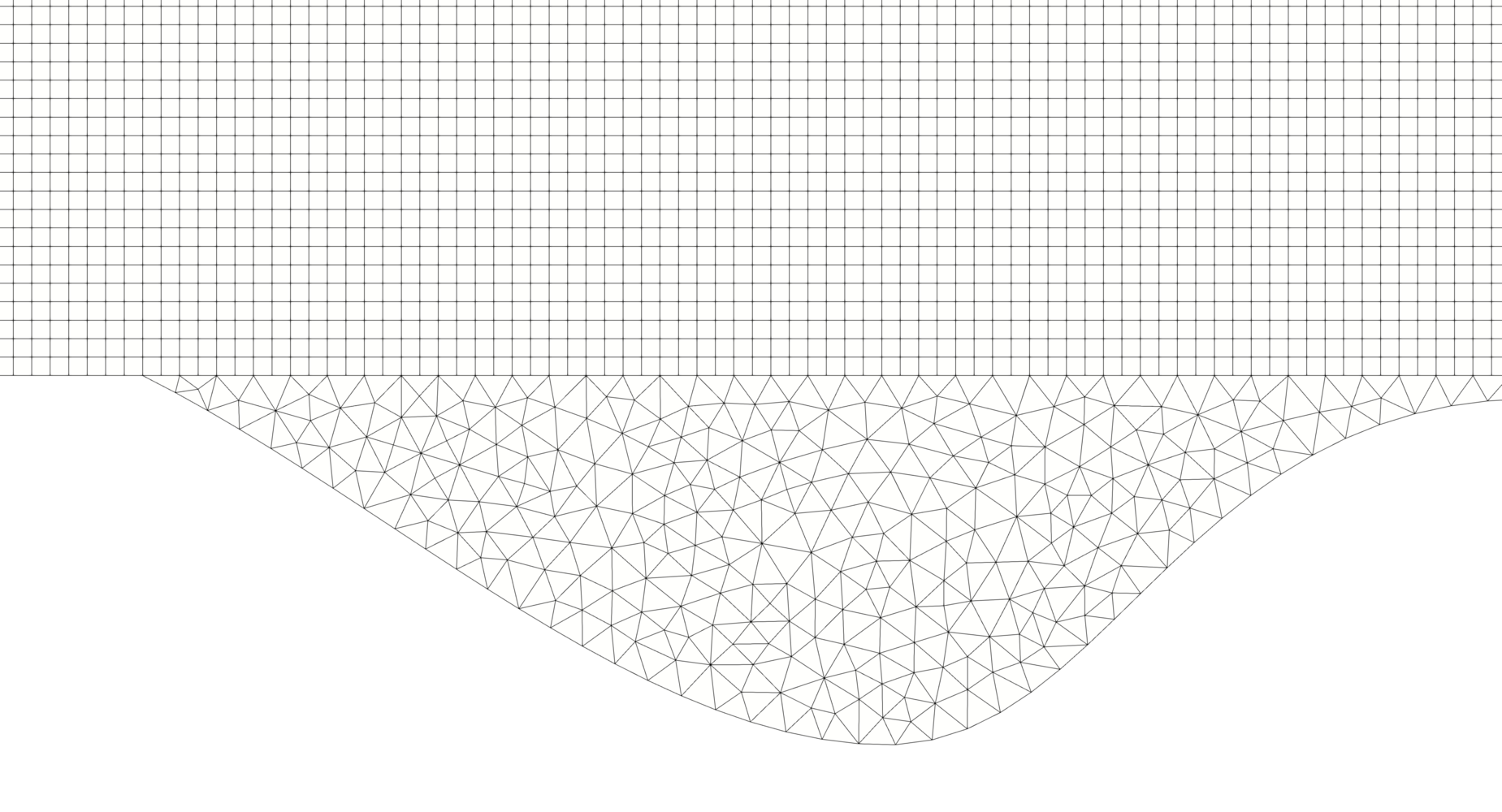}
	\caption{Detailed view on the grid used in test case three.}
	\label{pc:testTwoGrid}
\end{figure}

Another advantage of using MPFA in the porous medium domain is the ability to use unstructured grids while maintaining
 consistency of the scheme. The test case presented in this section considers a porous medium domain in which geometrical
 constraints favor the use of triangles (unstructured) over quadrilaterals for its discretization. This situation can arise,
 for example, in environmental applications considering the exchange processes between the atmosphere and the subsurface, where
 the latter can be composed of complex shaped geological layers. An illustration of the setup of this test case is given in
 Figure \ref{pc:settingCaseTwo}, while a detailed view on the discretization of the two compartments is depicted in Figure \ref{pc:testTwoGrid}.

\begin{figure}[ht!]
	\centering
	\includegraphics[width=0.85\linewidth,keepaspectratio]{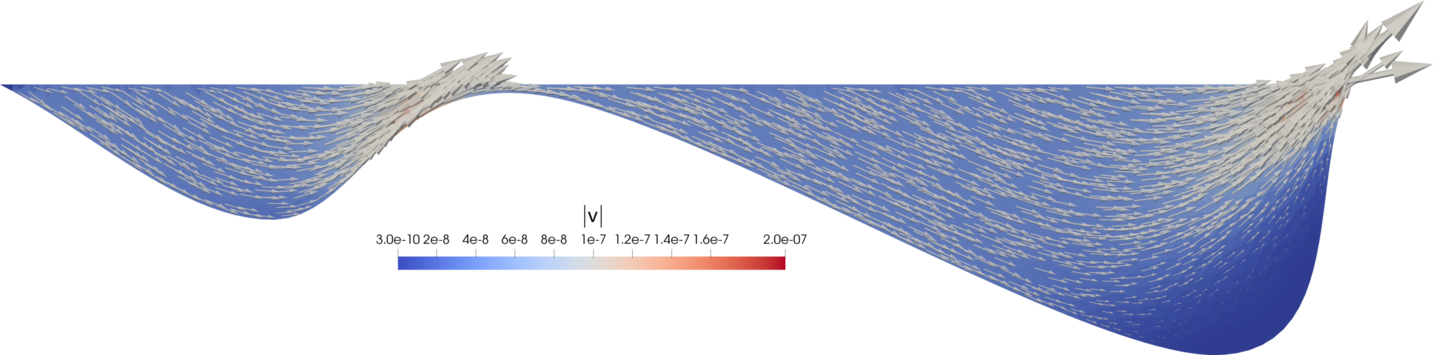}
	\includegraphics[width=0.85\linewidth,keepaspectratio]{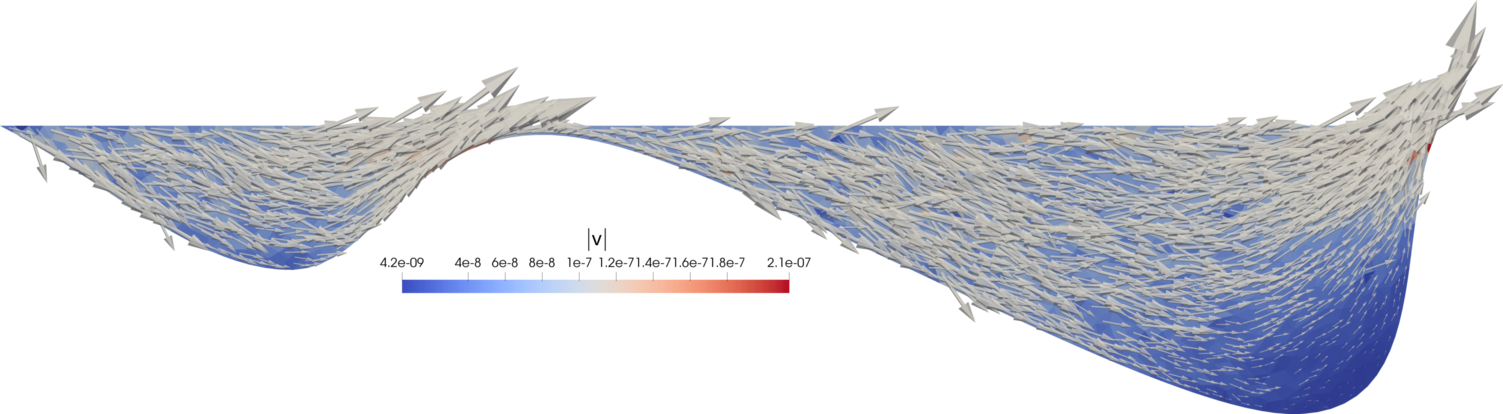}
	\caption{Velocity distribution in the porous medium for the third test case with $\beta = 10$ and $\alpha = -45^\circ$.
             The upper and lower pictures depict the result using MPFA and TPFA in the porous medium domain, respectively.}
	\label{pc:testTwoPMVelocity_-45}
\end{figure}

\begin{figure}[ht!]
	\centering
    \includegraphics[width=0.85\linewidth,keepaspectratio]{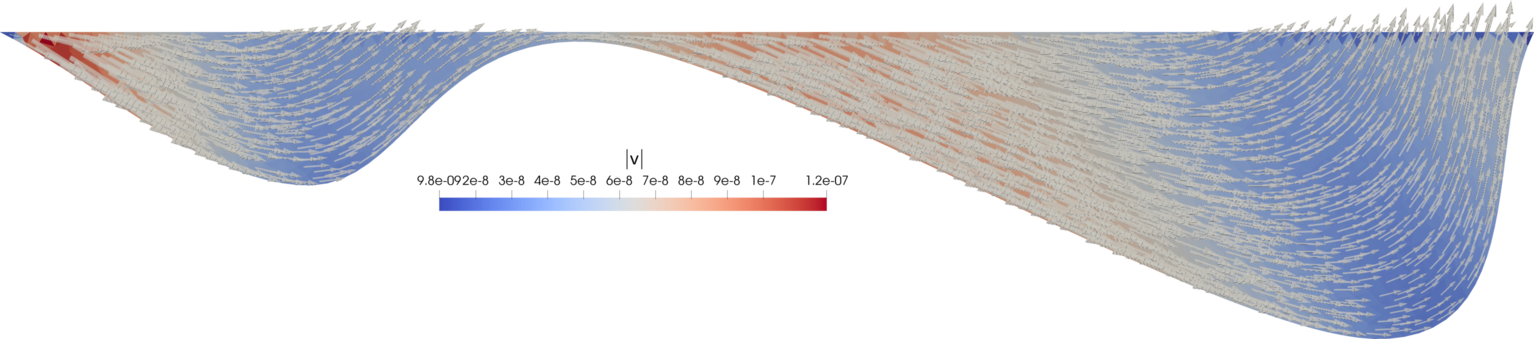}
	\includegraphics[width=0.85\linewidth,keepaspectratio]{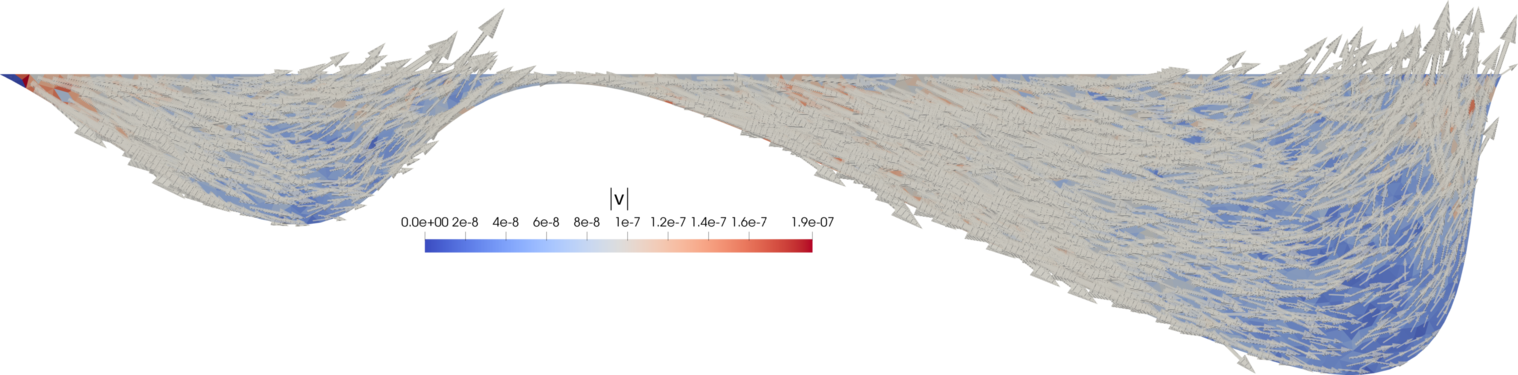}
	\caption{Velocity distribution in the porous medium for the third test case with $\beta = 10$ and $\alpha = 45^\circ$.
             The upper and lower pictures depict the result using MPFA and TPFA in the porous medium domain, respectively.}
	\label{pc:testTwoPMVelocity_45}
\end{figure}

Two sets of simulations were performed using the rotation angles $\alpha = 45^\circ$ and $\alpha = -45^\circ$, together with $\beta = 10$,
 $\unit[k=10^{-6}]{m^2}$ and $\Delta p = \unit[1 \cdot 10^{-6}]{Pa}$. The results using both MPFA and TPFA in the porous medium domain are shown
 in the Figures \ref{pc:testTwoPMVelocity_-45} and \ref{pc:testTwoPMVelocity_45}, respectively. Please note that unlike before, the velocity vectors are now scaled by magnitude.
 As expected, the solutions using TPFA exhibit
 a more distorted velocity field originating from the scheme being inconsistent on both unstructured grids and for anisotropic tensors.
 On the other hand, the solutions using MPFA provide velocity fields that follow the geometrical features of the porous medium.
 For $\alpha = -45^\circ$, two main regions in which the flux from the porous medium into the free-flow domain is concentrated,
 and where the maximum velocities occur, can be observed (see Figure \ref{pc:testTwoPMVelocity_-45}). The direction of these
 maximum velocities coincides with the direction of highest permeability. In contrast to that, the inflow from the free-flow domain into the porous
 medium occurs in the direction of the lowest permeability and thus at smaller velocities.
 In the case of $\alpha = 45^\circ$, the highest velocities are observed in the regions where an inflow from the free-flow domain into the porous
 medium occurs, again following the direction of the highest permeability (see Figure \ref{pc:testTwoPMVelocity_45}). With TPFA being used in the porous
 medium domain, the differences between the velocity fields obtained from the two angles turn out to be smaller in comparison to the solutions obtained with MPFA.
 Furthermore, the low velocity regions seem to be generally overestimated. This effect shows itself also in the integrated transfer flux across the interface, which
 results in $\unit[1.3 \cdot 10^{-8}]{kg/s}$ for MPFA and $\unit[2.2 \cdot 10^{-8}]{kg/s}$ for TPFA, thus, around $\unit[70]{\%}$ higher.

\section{Conclusions}
\label{sec:conclusions}

In this work, a discretization method is proposed for problems concerning free flow coupled to
porous medium flow. The method combines the stability of the staggered grid finite volume method
for the free-flow equations with the consistency of MPFA finite volume methods for flows in
anisotropic porous media. We have shown how appropriate alignment of the grids results in a natural
coupling between the two discretization schemes due to coinciding degrees of freedom.

The stability and consistency of the method are emphasized with the use of numerical experiments.
Especially in the presence of anisotropy in the porous medium, a significant difference is
observed with respect to the widely used, but inconsistent, TPFA finite volume method. We moreover
emphasized the use of unstructured grids in the porous medium, allowing for computations on more
general geometries. The use of unstructured grids in the free-flow domain \cite{eymard2014a} is a topic for future
investigation.

Future work will moreover address the extension of the presented model to multi-phase flow simulations including compositional and non-isothermal effects,
 with a special focus on evaporative processes at the interface between the porous medium and the free-flow domain. This extension
 would make the model applicable to a large variety of applications, as e.g.\ the drying of soil due to evaporation \cite{mosthaf2011a} and subsequent soil
 salinization \cite{jambhekar2015a}. This drying process can lead to the creation of fractures within the porous medium. Therefore, we want to employ
 a discrete fracture model for the description of the porous medium domain, for which there is a preexisting implementation available
 in \Dumux (see \cite{glaser2017discrete}). To increase efficiency, we will also use other finite volume schemes for the porous domain, as for example those that have been recently presented in \cite{Schneider.ea:2018}.

\section*{Acknowledgements}
We thank the German Research Foundation (DFG) for supporting this work by funding SFB 1313, Research Project A02.
\bibliography{literature.bib}

\end{document}